\documentclass[12pt,amstags,fleqn]{article}

\usepackage{url}

\usepackage{amssymb}

\usepackage{amsmath}
\usepackage{amsthm}

\usepackage{xspace}

\theoremstyle{plain}
\newtheorem{theorem}{Theorem}[section]

\theoremstyle{definition}

\newtheorem{construction}[theorem]{Construction}
\newtheorem{example}[theorem]{Example}
\newtheorem{remark}[theorem]{Remark}

\newcommand{\darts}{\Omega}

\def\ll{{\textstyle \ast}}
\def\rr{{\scriptscriptstyle \triangle}}

\newcommand{\opa}{\ll} 
\newcommand{\opb}{\rr}

\newcommand{\eq}{\textnormal{Eq}}

\usepackage{ifpdf}

\ifpdf
    \usepackage[pdftex]{graphicx}
\else
    \usepackage{graphicx}
\fi

\begin{document}

\title{An enumeration of equilateral triangle dissections}
\author{\Large 
Ale\v s Dr\'apal\footnote{Supported by grant MSM 0021620839}\\
Department of Mathematics \\ Charles University \\ Sokolovsk\'a 83 \\ 186 75 Praha 8 \\ Czech Republic \\
~ \\
Carlo H\"{a}m\"{a}l\"{a}inen\footnote{Supported by Eduard \v Cech center, grant LC505.}  \\
Department of Mathematics \\ Charles University \\ Sokolovsk\'a 83 \\ 186 75 Praha 8 \\ Czech Republic \\
{\texttt carlo.hamalainen@gmail.com}\\
}
%
%
%

\maketitle

\begin{abstract}
We enumerate all dissections of an equilateral triangle into smaller
equilateral triangles up to size $20$, where each triangle has integer
side lengths. A {\em perfect} dissection has no two triangles of the
same side, counting up- and down-oriented triangles as different.
We computationally prove W.~T.
Tutte's conjecture that the
smallest perfect dissection has size $15$
and we find all perfect
dissections up to size $20$.
\end{abstract}

\section{Introduction}

\begin{figure}[hbt]
\begin{center}
\includegraphics{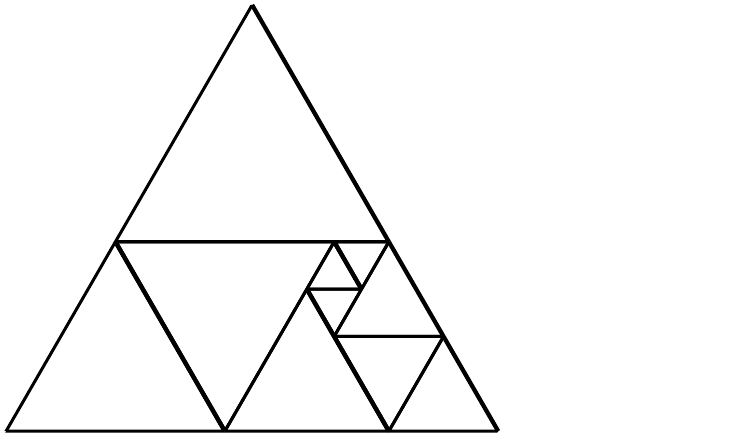}
\end{center}
\caption{An example of an equilateral triangle dissection.}
\label{exDissection}
\end{figure}

We are concerned with the following problem: given an equilateral
triangle $\Sigma$, find all dissections of $\Sigma$ into smaller
nonoverlapping equilateral triangles. The {\em size} of a dissection is
the number of nonoverlapping equilateral triangles. An example of such a dissection
of size $10$ is shown in Figure~\ref{exDissection}. 
It is well known that in such a dissection all triangles may be
regarded as triangles with sides of integer length.
Dissections of squares have been studied earlier~\cite{MR0003040}
as well as dissections of squares into right-angled isosceles
triangles~\cite{MR1794696}. Recently, 
Laczkovich\cite{MR1092545} studied tilings of polygons by similar
triangles. The earliest study of dissections of equilateral triangles
into equilateral triangles is by Tutte~\cite{MR0027521}. The problem of
dissecting a triangle is
different to normal tiling problems where the size of the tiles is known
in advance and the tiling area may be infinite.
A naive approach to enumerating dissections is to first fix the sizes and
number of
the dissecting triangles. Observe that in any dissection, some triangles
will be oriented in the same way as the triangle $\Sigma$ (these are the
up-triangles) while the oppositely oriented triangles are the down-triangles.
Let $u_{s}$ and $d_{s}$ be the number of up and down triangles of side
length $s$, respectively. For any down-triangle the horizontal side is
adjacent to the horizontal side of some number of up-triangles. The
up-triangles along the bottom of $\Sigma$ are not adjacent to any
down-triangle. So if the triangle $\Sigma$ has side length $n \in
\mathbb{N}$ then
\begin{equation}\label{eqnSideSum}
\sum_{s} s u_{s} = \sum_{s} s d_{s} + n.
\end{equation}
A triangle with length ${s}$ has height $\sqrt{3}{s}/2$ and so 
the triangle areas give the relation
\begin{equation}\label{eqnAreaSum}
\sum_s u_s s^2 + \sum_s d_s s^2 = n^2.
\end{equation}
For small values of $n$ we can solve \eqref{eqnSideSum} and \eqref{eqnAreaSum}
for the permissible
size and number of up and down triangles, and this data may guide an
exhaustive search.
We will consider up and down triangles not to be congruent even
if they are of the same size. A {\em perfect tiling} or
{\em perfect dissection} has no pair of congruent triangles.
This definition of a perfect dissection arises from the fact that 
it is impossible to have a perfect dissection if orientation is
ignored~\cite{MR0003040}. This fact can also easily derived from the
results in~\cite{aleshamming}.
Tutte conjectured~\cite{MR0003040,MR0027521} that the smallest perfect
dissection has size $15$ (see also~\cite{squaring}).  Unfortunately
solving \eqref{eqnSideSum} and \eqref{eqnAreaSum} with $n = 15$ is computationally intensive and so
another approach is needed.  Using our enumerative methods we confirm Tutte's
conjecture and provide all perfect dissections up to size $20$.

Lines parallel to the outer sides of the main triangle
that are induced by a side
of a dissecting triangle are \emph{dissecting lines}.
For any dissecting line $l$,
the union of all sides of dissecting triangles
that are incident to $l$ forms
one or more contiguous segments. If there are two or more segments,
then on $l$ there exist
two dissecting vertices such that all
triangles in between are cut by the line into two parts. If such
a situation arises for no dissecting line and if no
dissecting vertex is incident to six dissecting triangles,
then we call the dissection {\em separated}.
We enumerate all isomorphism classes of 
separated and nonseparated dissections up to size $20$.

\section{Dissections and latin bitrades}

The connection between equilateral triangle dissections and latin
bitrades was first studied in~\cite{aleshamming}. The presentation here
follows~\cite{alesdissections}.  
Consider an equilateral triangle $\Sigma$ that is dissected
into a finite number of equilateral triangles. 
Dissections will be always assumed to be nontrivial so the
number of dissecting triangles is at least four. 
Denote by $a$, $b$ and $c$
the lines induced by the sides of $\Sigma$.
Each side of a dissecting triangle has to be parallel
to one of $a$, $b$, or $c$. If $X$ is a vertex
of a dissecting triangle, then $X$ is a vertex of exactly one,
three or six dissecting triangles. Suppose that there is
no vertex with six triangles and consider triples 
$(u,v,w)$ of lines that are parallel to $a$, $b$ and $c$, respectively,
and meet in  a vertex of a dissecting triangle that is not a vertex
of $\Sigma$. The set of all these triples together with the triple
$(a,b,c)$ will be denoted by $T^\ll$, and by $T^\rr$ we shall
denote the set of all triples $(u,v,w)$ of lines that are yielded 
by sides of a dissecting triangle (where $u$, $v$ and $w$ are again parallel to
$a$, $b$ and $c$, respectively). The following conditions
hold:
\begin{enumerate}
\item[(R1)] Sets $T^\ll$ and $T^\rr$ are disjoint;
\item[(R2)] for all $(p_1,p_2,p_3)\in T^\ll$ and 
all $r,s \in \{1,2,3\}$, $r \ne s$, there exists exactly one
$(q_1,q_2,q_3) \in T^\rr$ with $p_r = q_r$ and $p_s = q_s$; and
\item[(R3)] for all $(q_1,q_2,q_3)\in T^\rr$ and
all $r,s \in \{1,2,3\}$, $r \ne s$, there exists exactly one
$(p_1,p_2,p_3) \in T^\ll$ with $q_r = p_r$ and $q_s = p_s$.
\end{enumerate}

Note that (R2) would not be true if there had existed six dissecting
triangles with a common vertex. Conditions (R1--3) are, in fact, axioms
of a combinatorial object called latin bitrades
\cite[p.~148]{wanlesshandbook}. A bitrade is usually denoted
$(T^{\opa},\, T^{\opb})$. Observe that the bitrade $(T^{\opa},\,
T^{\opb})$ associated with a dissection encodes
qualitative
(structural) information about the segments and intersections of
segments in the dissection. The sizes and number of dissecting triangles
can be recovered by solving a system of equations derived from 
the bitrade (see below).

Dissections are related to a class of latin bitrades with genus~$0$. To
calculate the genus of a bitrade we use a permutation representation to
construct an oriented combinatorial surface~\cite{Dr9,hamalainen2007,ales-geometrical}.
For $r \in \{1,\, 2,\, 3\}$, define the map 
$\beta_r \colon T^{\opb} \rightarrow
T^{\opa}$ where $(a_1,\, a_2,\, a_3) \beta_r = (b_1,\, b_2,\, b_3)$ if and only if
$a_r \neq b_r$
and $a_i = b_i$ for $i \neq r$.
By (R1-3)
each $\beta_r$ is a bijection.
Then
$\tau_1,\, \tau_2,\, \tau_3\colon T^{\opa} \rightarrow
T^{\opa}$ are defined by
\begin{align}
\tau_1 &= \beta_2^{-1}\beta_3, \qquad
\tau_2 = \beta_3^{-1}\beta_1, \qquad
\tau_3 = \beta_1^{-1}\beta_2. \label{eqnTau}
\end{align}
We refer to
$[\tau_1,\, \tau_2,\, \tau_3]$
as the $\tau_i$ {\em representation}. To get a combinatorial surface
from a bitrade we use the following construction:

\begin{construction}\label{constructionTauiSurface}
Let $[\tau_1,\, \tau_2,\, \tau_3]$ be the representation for 
a bitrade where the $\tau_i$ act on
the set $\darts$. Define
vertex, directed edge, and face sets by:
\begin{alignat*}{2}
V &= \darts \\
E &= \{ (x,y) \mid x \tau_1 = y \} \cup \{ (x,y) \mid x \tau_2 = y \} \cup \{ (x,y) \mid x \tau_3 = y \} \\
F &= \{ (x,y,z) \mid \text{ $x \tau_1 = y$, $y \tau_2 = z$, $z \tau_3 = x$} \} \\
  & \qquad \cup \{ (x_1,\, x_2,\, \dots,\, x_r) \mid \text{ $(x_1,\,
x_2,\, \dots,\, x_r)$ is a cycle of $\tau_1$ } \} \\
  & \qquad \cup \{ (x_1,\, x_2,\, \dots,\, x_r) \mid \text{ $(x_1,\,
x_2,\, \dots,\, x_r)$ is a cycle of $\tau_2$ } \} \\
  & \qquad \cup \{ (x_1,\, x_2,\, \dots,\, x_r) \mid \text{ $(x_1,\,
x_2,\, \dots,\, x_r)$ is a cycle of $\tau_3$ } \}
\end{alignat*}
where $(x_1,\, x_2,\, \dots,\, x_k)$ denotes a face with $k$~directed edges
$(x_1,\, x_2)$, $(x_2,\, x_3)$, $\dots$,
$(x_{k-1},\, x_k)$, $(x_k,\, x_1)$ for vertices
$x_1,\, \dots,\, x_k$.

The first set in the definition of $F$ is the set of {\em triangular
faces}, while the other three are the {\em $\tau_i$ faces}.
Assign triangular faces a positive (anticlockwise) orientation, and
assign
$\tau_i$ faces negative (clockwise) orientation. Now glue the faces
together where they share a common directed edge $x \tau_i = y$,
ensuring that edges come together with opposite orientation.
\end{construction}

\begin{figure}[htb]
\begin{center}
\includegraphics{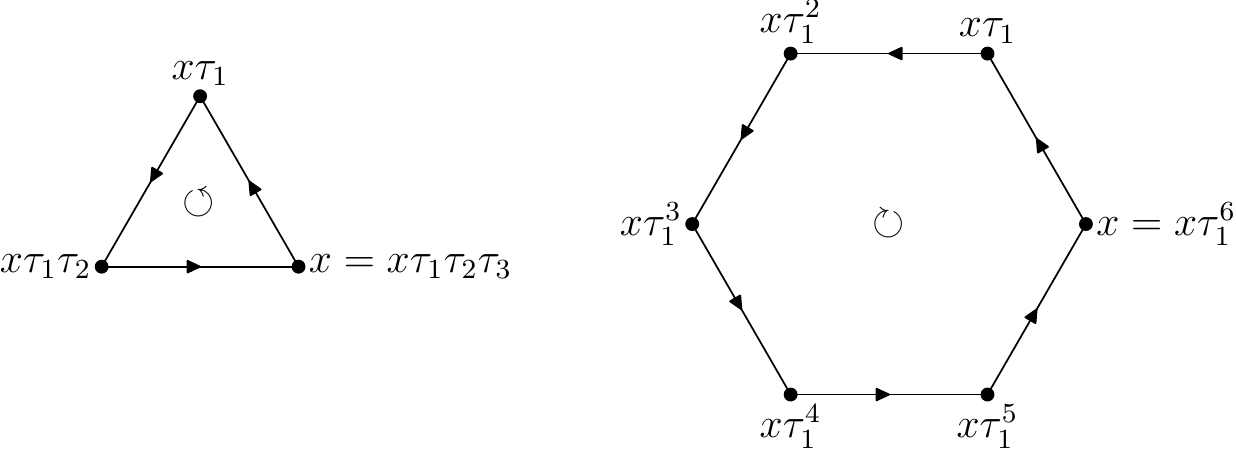}
\end{center}
\caption{Orienting faces in the combinatorial surface.}
\label{figTwoOrientations}
\end{figure}

\begin{figure}[htb]
\begin{center}
\includegraphics{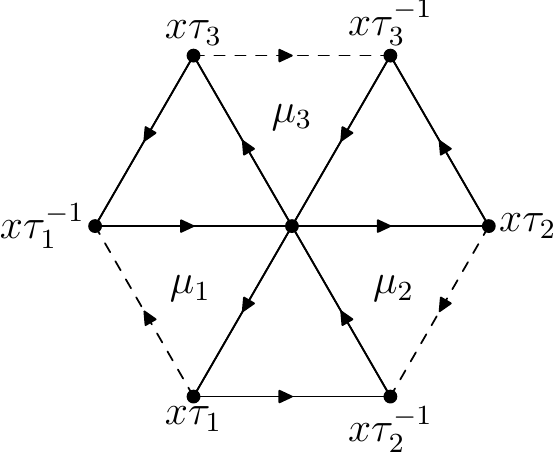}
\end{center}
\caption{Rotation scheme for the combinatorial surface. Dashed
lines represent one or more edges, where $\mu_i$ is a cycle of $\tau_i$.
The vertex in the centre is $x$.}
\label{figDrapalRotationScheme}
\end{figure}

The orientation of a triangular face is
shown in Figure~\ref{figTwoOrientations} along with the orientation for
a $\tau_i$ face. For the sake of concreteness we have illustrated a
$6$-cycle face due to a $6$-cycle of $\tau_1$. 
Figure~\ref{figDrapalRotationScheme} shows
the rotation scheme for an arbitrary vertex in the surface.

For a bitrade $T = (T^{\opa},\, T^{\opb})$ with representation 
$[\tau_1,\, \tau_2,\, \tau_3]$
on the set $\darts$,
define
$\text{order}(T)  = z(\tau_1) + z(\tau_2) + z(\tau_3)$, the total number
of cycles, and
$\text{size}(T) = \left| \darts \right|$, the total number of points 
that the $\tau_i$ act on. 
By some basic counting arguments we find that there are 
$\text{size}(T)$ vertices, $3 \cdot \text{size}(T)$ edges, and 
$\text{order}(T) + \text{size}(T)$ faces. Then Euler's formula
$V-E+F=2-2g$ gives
\begin{equation}\label{eqnDrapalGenus}
\text{order}(T) = \text{size}(T) + 2 - 2g
\end{equation}
where $g$ is the genus of the combinatorial surface. We say that
the bitrade $T$ has genus $g$.
A {\em spherical bitrade} has genus~$0$.

\begin{example}\label{exspherical}
The following bitrade is spherical:
\begin{align*}
T^{\opa} = 
\begin{array}{|c||c|c|c|c|c|}
\hline
\opa & 0 & 1 & 2 & 3 & 4\\
\hline
 \hline 0 & 0 & ~ & 2 & ~ & 4\\
 \hline 1 & ~ & ~ & ~ & 4 & 2\\
 \hline 2 & 1 & 3 & 0 & 2 & ~\\
 \hline 3 & 4 & 1 & ~ & 3 & ~\\
\hline
\end{array}
& \quad & 
T^{\opb} = 
\begin{array}{|c||c|c|c|c|c|}
\hline
\opb & 0 & 1 & 2 & 3 & 4\\
\hline
 \hline 0 & 4 & ~ & 0 & ~ & 2\\
 \hline 1 & ~ & ~ & ~ & 2 & 4\\
 \hline 2 & 0 & 1 & 2 & 3 & ~\\
 \hline 3 & 1 & 3 & ~ & 4 & ~\\ 
\hline
\end{array}
\end{align*}
Here, the $\tau_i$ representation is
\begin{align*}
\tau_1 &= (000, 022, 044)(134, 142) (201,213,232,220) (304,333,311) \\
\tau_2 &= (000,304,201)(213,311)(022,220)(134,232,333)(044,142) \\
\tau_3 &= (000,220)(201,311)(022,232,142)(213,333)(044,134,304)
\end{align*}
where the triples $ijk$ refer to entries $(i,j,k) \in T^{\opa}$.
\end{example}

We will generally assume that a bitrade is {\em separated}, that is,
each row, column, and symbol is in bijection with a cycle of $\tau_1$,
$\tau_2$, and $\tau_3$, respectively.

We now describe how to go from a separated spherical latin bitrade to a
triangle dissection (for more details see~\cite{alesdissections}).  
Let $T= (T^\ll, T^\rr)$ be a latin bitrade. It is natural to have
different unknowns for rows, columns and symbols, and so we assume
that $a_i \ne b_j$ whenever $(a_1,a_2,a_3)$, $(b_1,b_2,b_3) \in T^\ll$
and $1 \le i < j \le 3$.  (If the condition is violated,  then $T$
can be replaced by an isotopic bitrade for which it is satisfied.)
Fix a triple $a = (a_1,a_2,a_3) \in T^\ll$ and form the set of equations
$\eq(T)$ consisting of
$a_1=0$, $a_2=0$, $a_3 =1$
and 
$b_1+b_2 = b_3$  if $(b_1,b_2,b_3) \neq (a_1,a_2,a_3)$
and $(b_1,b_2,b_3) \in T^\ll$.
The theorem below shows that if $T$ is a spherical latin
bitrade then $\eq(T, a)$ always has a unique solution in the rationals.
The pair $(T,a)$ will be called a \emph{pointed} bitrade.

Write
$\bar r_i$,
$\bar c_j$,
$\bar s_k$
for the solutions in $\eq(T, a)$ for row variable $r_i$, column variable
$c_j$, and symbol variable $s_k$, respectively. We say that 
a solution to $\eq(T, a)$ is {\em separated} if 
$\bar r_i \neq \bar r_{i'}$ whenever $i \neq i'$ (and similar for
columns and symbols).

For each entry $c = (c_1,c_2,c_3) \in T^\rr$ we form the triangle
$\Delta(c,a)$ which is bounded by the lines
$y = \bar c_1$,
$x = \bar c_2$,
$x + y = \bar c_3$.
Of course, it is not clear that $\Delta$ is really
a triangle, i.~e.~that the three lines do not meet in a single
point. If this happens, then we shall say that $\Delta(c,a)$ {\em degenerates}.
Let $\Delta(T,a)$ denote the subset of $T^\rr$ such that
$\Delta(c,a)$ does not degenerate.

A separated dissection 
with $m$ vertices corresponds to a separated spherical bitrade
$(T^{\opa},\, T^{\opb})$ where $\left| T^{\opa} \right| = m-2$.  
One of the main results of~\cite{alesdissections} is the following
theorem:

\begin{theorem}[\cite{alesdissections}]\label{sepdissections}
Let $T=(T^\ll, T^\rr)$ be a spherical latin bitrade,
and suppose that $a=(a_1,a_2,a_3) \in T^\ll$ is a triple such
that the solution to $\eq(T,a)$ is separated. Then the set
of all triangles $\Delta(c,a)$, $c \in T^\rr$, dissects the triangle
$\Sigma = \{(x,y);$ $x\ge 0$, $y\ge 0$ and $x+y \le 1\}$. This
dissection is separated.
\end{theorem}

An equilateral dissection can be obtained by applying the transformation
$(x,y) \mapsto (y/2 + x, \sqrt 3 y/2)$.

\begin{figure}[htbp]
\begin{center}
\includegraphics[scale=0.75]{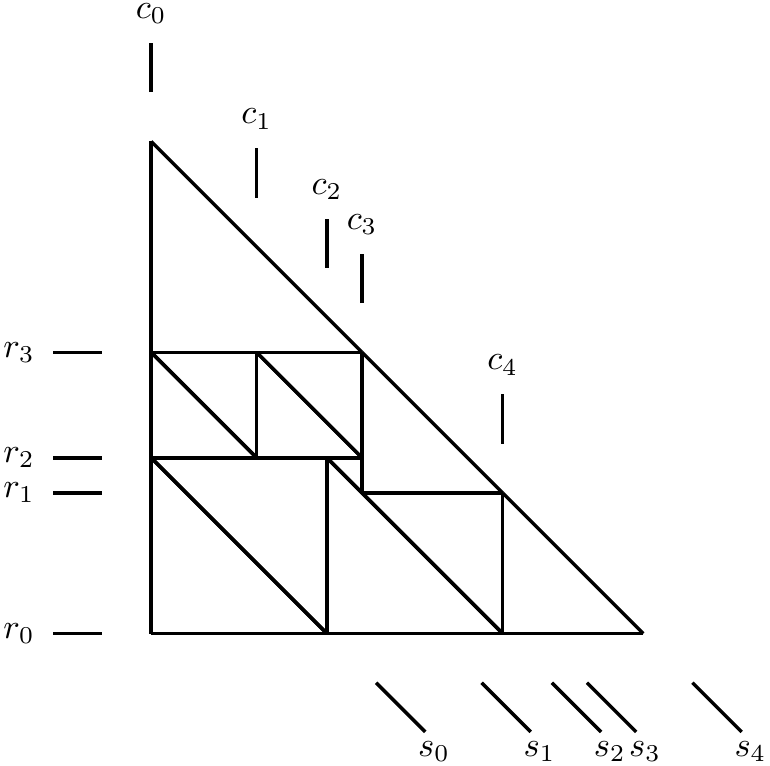}
\end{center}
\caption{Separated dissection for a spherical bitrade. The labels $r_i$, $c_j$,
$s_k$, refer to lines $y = \bar r_i$,
$x = \bar c_j$,
$x+y = \bar s_k$, respectively. The trade $T^{\opa}$ has $12$ entries
and the dissection has $12 + 2 = 14$ vertices.
Applying the transformation $(x,y) \mapsto (y/2 + x, \sqrt 3 y/2)$ gives
an equilateral dissection.}
\label{fg2}
\end{figure}

\begin{example}
Consider the following spherical bitrade 
$(T^{\ll},\, T^{\rr})$:  \\ 
\[
T^{\ll} = 
\begin{array}{|c||c|c|c|c|c|}
\hline \ll & c_0 & c_1 & c_2 & c_3 & c_4\\
\hline \hline r_0 & s_4 & ~ & s_0 & ~ & s_2\\
\hline r_1 & ~ & ~ & ~ & s_2 & s_4\\
\hline r_2 & s_0 & s_1 & s_2 & s_3 & ~\\
\hline r_3 & s_1 & s_3 & ~ & s_4 & ~\\
\hline \end{array} 
\quad 
T^{\rr} = 
\begin{array}{|c||c|c|c|c|c|}
\hline \rr & c_0 & c_1 & c_2 & c_3 & c_4\\
\hline \hline r_0 & s_0 & ~ & s_2 & ~ & s_4\\
\hline r_1 & ~ & ~ & ~ & s_4 & s_2\\
\hline r_2 & s_1 & s_3 & s_0 & s_2 & ~\\
\hline r_3 & s_4 & s_1 & ~ & s_3 & ~\\
\hline \end{array}
\label{eqnsphericalbitrade}
\]

Let $a = (a_1,a_2,a_3) = (r_0,c_0,s_4)$. Then the system of equations
$\textnormal{Eq}(T, a)$ has the solution
\begin{align*}
\bar r_0 &= 0,\, \bar r_1 = 2/7,\, \bar r_2 = 5/14,\, \bar r_3 = 4/7 \\
\bar c_0 &= 0,\, \bar c_1 = 3/14,\, \bar c_2 = 5/14,\, \bar c_3 = 3/7,\, \bar c_4 = 5/7 \\
\bar s_0 &= 5/14,\, \bar s_1 = 4/7,\, \bar s_2 = 5/7,\, \bar s_3 = 11/14,\, \bar s_4 = 1.
\end{align*}
The dissection is shown in Figure~\ref{fg2}.
 Entries of
$T^{\rr}$ correspond to triangles in the dissection. For example,
$(r_0,c_0,s_0) \in T^{\rr}$ is the triangle bounded by the lines
$y = \bar r_0 = 0$, 
$x = \bar c_0 = 0$, 
$x+y = \bar s_0 = 5/14$ while 
$(r_1,c_3,s_2) \in T^{\ll}$ corresponds to the intersection of the
lines
$y = \bar r_1 = 2/7$, 
$x = \bar c_3 = 3/7$, 
$x+y = \bar s_2 = 5/7$.
\end{example}

\begin{figure}[hbt]
\begin{center}
\includegraphics{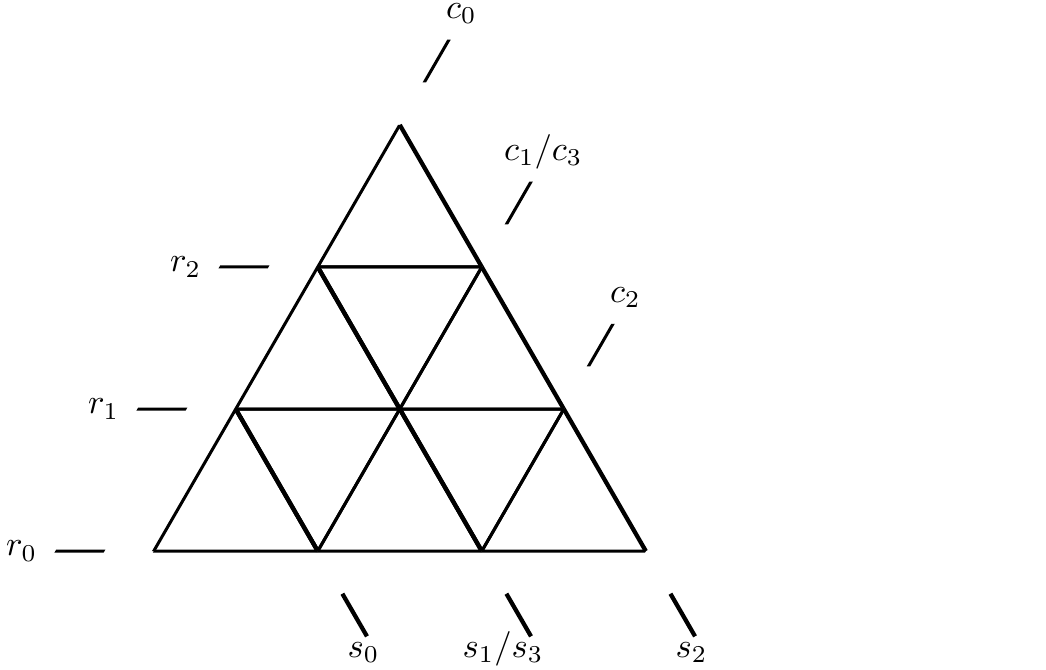}
\end{center}
\caption{A dissection with a vertex of degree $6$.}
\label{figsixway}
\end{figure}

\begin{remark}
Suppose that the dissection $\Sigma$ has no vertex of degree $6$.
Pick a vertex $X$ of degree~$4$. If we
move to the right along the row segment to the next vertex $X'$ then we
have $X \tau_1 = X'$. Similarly, moving along the diagonal segments
gives the action of $\tau_2$ and $\tau_3$. If we identify the three vertices of degree $2$
then the dissection encodes, geometrically,
the permutation representation of the bitrade.

If a dissection has a vertex of degree~$6$ then the dissecting triangles
do not (uniquely) define a partial latin square and hence do not encode
a latin bitrade. However, we can recover a separated
bitrade by the following procedure. For each vertex $X$ of degree~$6$,
choose one segment (say, the $r_i$ segment) to stay fixed. Then for the $c_j$
and $s_k$ segments, label the column segment below $X$ with a new name $c'$
and label the symbol segment below $X$ with a new name $s'$. For example,
the centre vertex in Figure~\ref{figsixway} results in the new labels
$c_3$ and $s_3$. The resulting separated bitrade is:
\[ T^{\ast},\, T^{\triangle} =\begin{array}{|c||c|c|c|c|}
\hline \ast & c_0 & c_1 & c_2 & c_3\\
\hline \hline r_0 & s_2 & ~ & s_3 & s_0\\\hline r_1 & s_0 & s_1 & s_2 & s_3\\\hline r_2 & s_1 & s_2 & ~ & ~\\\hline\end{array}\phantom{x}\begin{array}{|c||c|c|c|c|}
\hline \triangle & c_0 & c_1 & c_2 & c_3\\
\hline \hline r_0 & s_0 & ~ & s_2 & s_3\\\hline r_1 & s_1 & s_2 & s_3 & s_0\\\hline r_2 & s_2 & s_1 & ~ & ~\\\hline\end{array}\]
This procedure works for any number of vertices of degree~$6$, as long
as care is taken to only relabel column or symbol segments below a
vertex of degree~$6$ and not to relabel a segment more than
once.\footnote{For a concrete implementation, see 
\texttt{generate\_bitrade\_via\_geometric\_data} in\\
\texttt{triangle\_dissections.py} in \cite{dissections}}
\end{remark}

Recently Theorem~\ref{sepdissections} has been strengthened to cover
nonseparated as well as separated dissections:

\begin{theorem}[\cite{alesnew}]\label{strongertheorem}
Let $T=(T^\ll, T^\rr)$ be a spherical latin bitrade. Then for 
any $a=(a_1,a_2,a_3) \in T^\ll$, the set $\Delta(T,a)$ of non-degenerate triangles
dissects the triangle
$\Sigma = \{(x,y);$ $x\ge 0$, $y\ge 0$ and $x+y \le 1\}$. The
dissection may not be separated.
\end{theorem}

A fundamental consequence of the theorem is that any dissection of size $s$
can be derived from a pointed spherical bitrade $(T, a)$ of size $s$.
Note that the systems of equations $\eq(T, a)$ have also other
applications: both \cite{nickWanlessBlah} and \cite{alesdissections}
use them to show that every spherical latin trade can be embedded
into a finite abelian group.

\section{Computational results}

Cavenagh and Lison\v ek~\cite{planareulerian} showed that spherical
bitrades are equivalent to planar Eulerian triangulations. To enumerate
triangle dissections we
use plantri~\cite{plantri,plantri2} to enumerate all planar
Eulerian triangulations up to size~$20$ (we also note that
in~\cite{wanlessenumeration} all trades and bitrades have been
enumerated up size~$19$). We wrote a 
plugin~\cite{code} to output the equivalent spherical latin bitrade
$(U^{\opa}, U^{\opb})$ for each triangulation. For each such
$(U^{\opa},U^{\opb})$ we find all
solutions $\eq(T, a)$ for all $a \in T^{\opa}$ and compute each
dissection (in practise this is a list of triangles $\Delta(c,a)$ for
each $c = (c_1,c_2,c_3) \in T^{\rr}$).
To filter out isomorphic dissections we apply all six elements of
the symmetry group for a unit-side equilateral triangle (identity, two rotations,
and three reflections). The {\em canonical signature} of a dissection 
$\Delta(T,a)$ is
the ordered list $[(x,y) \mid (x,y) \textnormal{ is a
vertex of $\Delta(T,a)$}]$. We repeat the whole process for the bitrade $(U^{\opb},U^{\opa})$
because there are spherical latin bitrades where
$(U^{\opa},U^{\opb})$ is not isomorphic to $(U^{\opb},U^{\opa})$. The
final counts for the number of dissections up to isomorphism are found
by simply removing duplicate signatures.

While the solutions to $\eq(T, a)$ exist in the rationals, we find it
easier to work with the final equilateral dissections instead. We use the SymPy package~\cite{sympy} to perform
exact symbolic arithmetic on the canonical form of each dissection.

\subsection{Dissections and automorphism groups}

Using Theorems~\ref{sepdissections} and \ref{strongertheorem}
we have enumerated the number of isomorphism classes of
dissections of size $n \leq 20$. We also record $A(n,k)$, the
number of dissections of size $n$ with
automorphism group of order $k$. See Figures~\ref{figCountSeparated} and
\ref{figCountSeparatedAndNonseparated} for the data.

The referee raised the question of asymptotic behaviour. Let us denote
by $d_n$ the number of all dissections of size $n$.  Thus $d_4 = 1$,
$\dots$, $d_{13} = 574$, $\dots$, $d_{20} = 2674753$.  There are some
reasons to believe that $d_n$ can be estimated as $\sigma(n)^n$,
where $\sigma(n)$ is a slowly growing function.  The asymptotic
behaviour of $\sigma(n)$ is not clear yet and is a subject of ongoing
research. Note however that $d_n \ge e_n = (3.43)^{n-8}$ for every
$n\notin \{18,19\}$ such that $8 \le n \le 20$.  If we put $\mu_n =
e_n/d_n$, then in this interval the approximate values of $\mu_n$
are $0.33$, $0.38$, $0.51$, $0.65$, $0.74$, $0.83$, $0.89$, $0.94$,
$0.97$, $0.99$, $1.00$, $1.00$ and $0.99$. The fact that $e_{18}$
and $e_{19}$ are slightly greater than $d_{18}$ and $d_{19}$, while
$e_{20}$ is smaller than $d_{20}$, seems to be surprising.
\begin{figure}[htbp]
\begin{center}
\begin{tabular}{|c|c|c|c|c|c|}
\hline $n$ & \# dissections & $A(n,1)$ & $A(n,2)$ & $A(n,3)$ & $A(n,6)$ \\
\hline
\hline 4  & 1 & 0 & 0 & 0 & 1 \\
\hline 6  & 1 & 0 & 1 & 0 & 0 \\
\hline 7  & 2 & 0 & 1 & 0 & 1 \\
\hline 8  & 3 & 2 & 1 & 0 & 0 \\
\hline 9  & 8 & 4 & 4 & 0 & 0 \\
\hline 10 & 20 & 15 & 4 & 0 & 1 \\
\hline 11 & 55 & 47 & 8 & 0 & 0 \\
\hline 12 & 161 & 146 & 15 & 0 & 0 \\
\hline 13 & 478 & 460 & 17 & 0 & 1 \\
\hline 14 & 1496 & 1459 & 37 & 0 & 0 \\
\hline 15 & 4804 & 4746 & 58 & 0 & 0 \\
\hline 16 & 15589 & 15506 & 82 & 0 & 1 \\
\hline 17 & 51377 & 51223 & 154 & 0 & 0 \\
\hline 18 & 172162 & 171923 & 239 & 0 & 0 \\
\hline 19 & 583810 & 583426 & 383 & 0 & 1 \\
\hline 20 & 1998407 & 1997752 & 655 & 0 & 0 \\
\hline
\end{tabular}
\end{center}
\caption{Number of separated dissections of size $n$, up to isomorphism.
For each $n$, the column $A(n,k)$
records the number of dissections of size $n$ with
automorphism group of order $k$.}
\label{figCountSeparated}
\end{figure}
\begin{figure}[htbp]
\begin{center}
\begin{tabular}{|c|c|c|c|c|c|}
\hline $n$ & \# dissections & $A(n,1)$ & $A(n,2)$ & $A(n,3)$ & $A(n,6)$ \\
\hline
\hline 4  & 1 & 0 & 0 & 0 & 1 \\
\hline 6  & 1 & 0 & 1 & 0 & 0 \\
\hline 7  & 2 & 0 & 1 & 0 & 1 \\
\hline 8  & 3 & 2 & 1 & 0 & 0 \\
\hline 9  & 9 & 4 & 4 & 0 & 1 \\
\hline 10 & 23 & 15 & 7 & 0 & 1 \\
\hline 11 & 62 & 51 & 11 & 0 & 0 \\
\hline 12 & 188 & 162 & 25 & 0 & 1 \\
\hline 13 & 574 & 532 & 39 & 0 & 3 \\
\hline 14 & 1826 & 1745 & 81 & 0 & 0 \\
\hline 15 & 5953 & 5795 & 157 & 0 & 1 \\
\hline 16 & 19664 & 19380 & 277 & 2 & 5 \\
\hline 17 & 66049 & 65489 & 560 & 0 & 0 \\
\hline 18 & 224700 & 223625 & 1070 & 0 & 5 \\
\hline 19 & 771859 & 769851 & 1992 & 8 & 8 \\
\hline 20 & 2674753 & 2670755 & 3998 & 0 & 0 \\
\hline
\end{tabular}
\end{center}
\caption{Number of separated and nonseparated dissections of size $n$, up to isomorphism.
For each $n$, the column $A(n,k)$
records the number of dissections of size $n$ with
automorphism group of order $k$.}
\label{figCountSeparatedAndNonseparated}
\end{figure}

\subsection{Perfect dissections}

\begin{figure}[htbp]
\begin{center}
\includegraphics{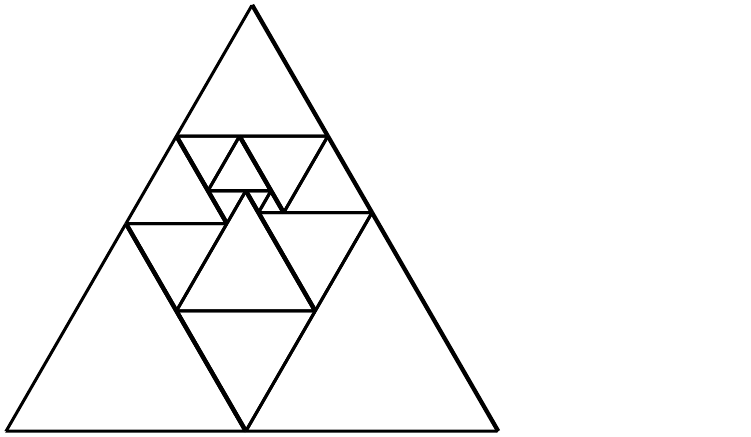}\includegraphics{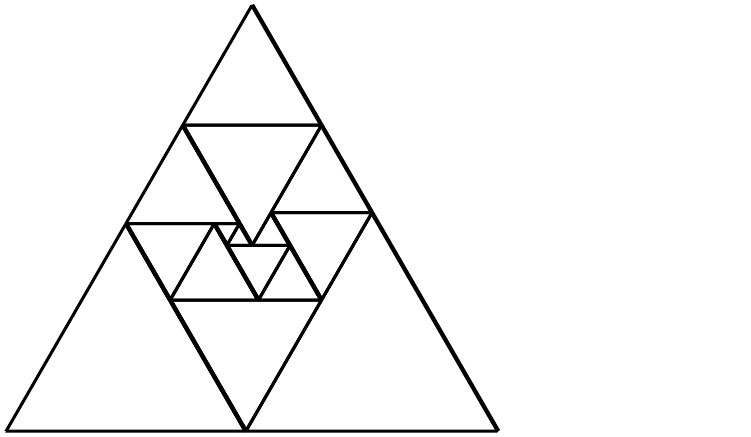}
\end{center}
\caption{The two perfect dissections of size $15$.}
\label{figPerfect15}
\end{figure}
\begin{figure}[htbp]
\begin{center}
\includegraphics{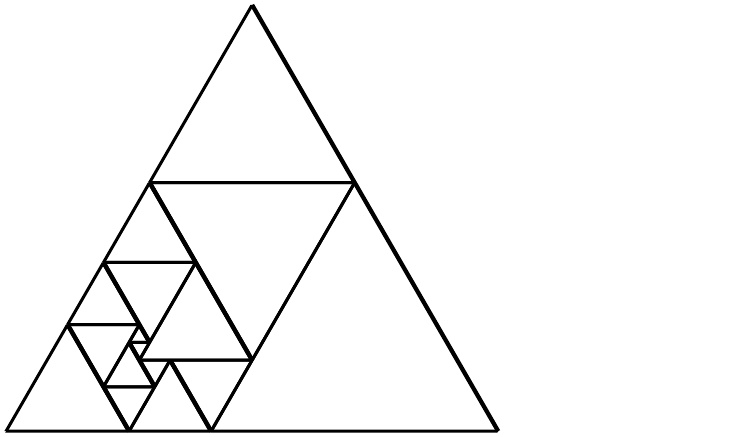}\includegraphics{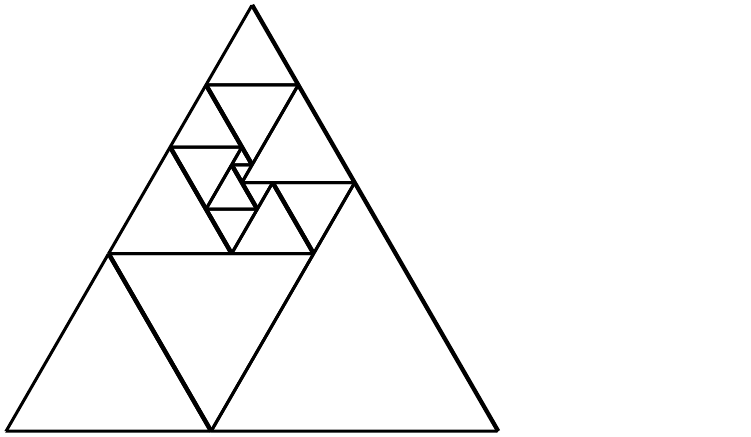} 
\end{center}
\caption{The two perfect dissections of size $16$.}
\label{figPerfect16}
\end{figure}
\begin{figure}[htbp]
\begin{center}
\includegraphics{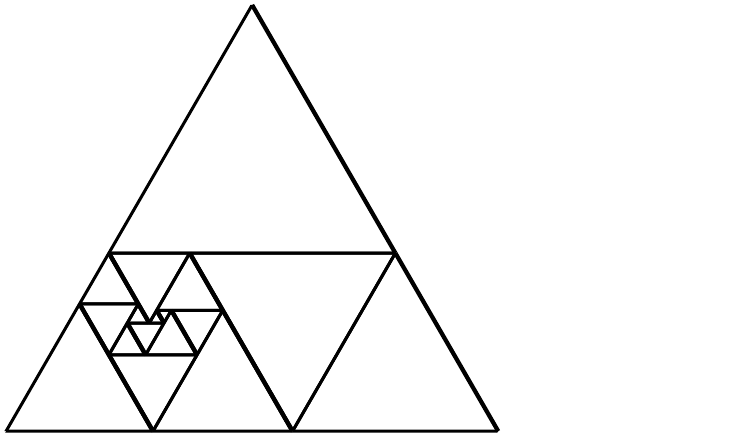}\includegraphics{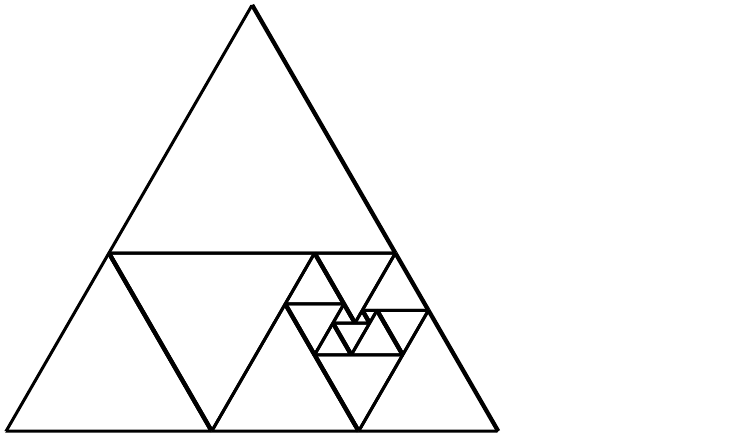}
\includegraphics{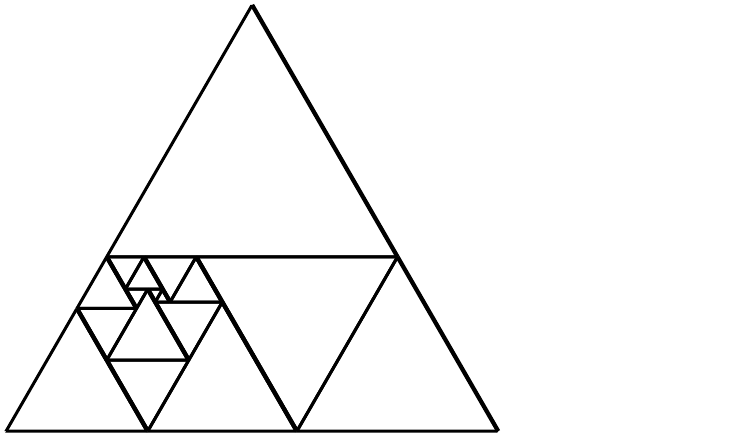}\includegraphics{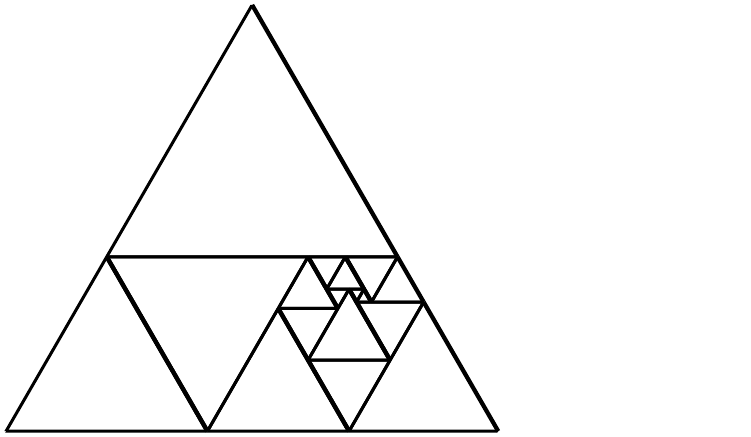}
\includegraphics{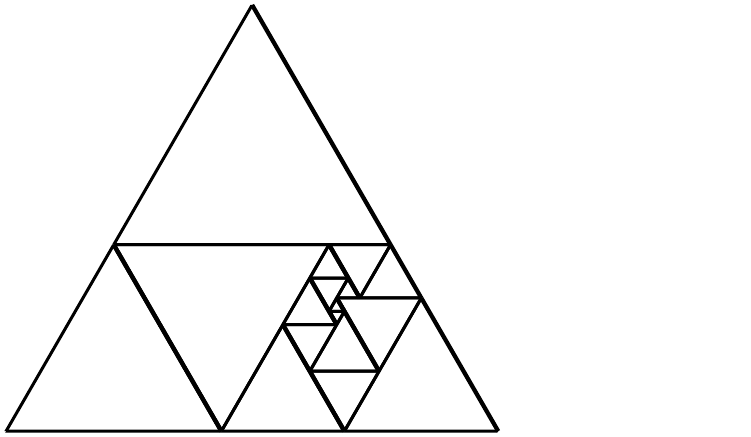}\includegraphics{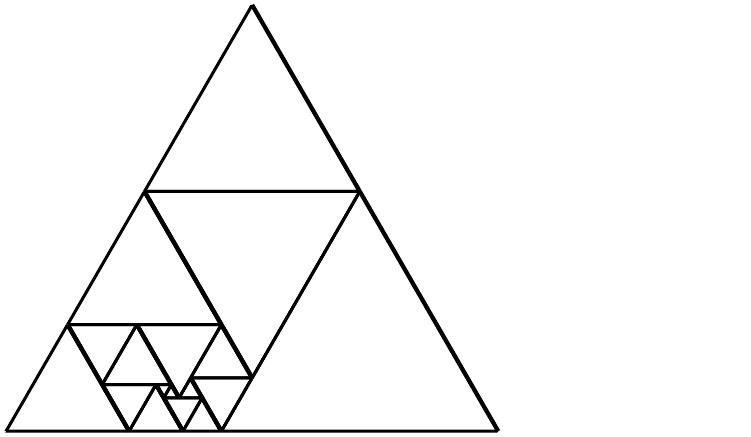}
\end{center}
\caption{The six perfect dissections of size $17$.}
\label{figPerfect17}
\end{figure}

Using our enumeration code we can confirm W.~T. Tutte's
conjecture~\cite{MR0003040,MR0027521}
that the smallest perfect dissection has size $15$ (see
Figure~\ref{figPerfect15}).
The perfect dissections of size $16$ and $17$ are shown in
Figures~\ref{figPerfect16} and \ref{figPerfect17}.
The perfect dissections of size up to $20$ are available in PDF format~\cite{dissections}.
The following table summarises the known number of isomorphism classes
of perfect dissections:

\begin{center}
\begin{tabular}{|c|c|}
\hline $n$ & \# perfect dissections \\
\hline
\hline 15 & 2 \\
\hline 16 & 2 \\
\hline 17 & 6 \\
\hline 18 & 23 \\
\hline 19 & 64 \\
\hline 20 & 181 \\
\hline
\end{tabular}
\end{center}

It is an open problem to determine if there exists a nonseparated
perfect dissection. If such a dissection exists then it will have size
greater than $20$.

\subsection{Other observations}

There exist bitrades of size $n \geq 10$ for which there is no
separated dissection of size $n$. The following table shows a sample
of these bitrades, giving the sizes of all possible dissections for
the particular bitrade of size $n$.

\begin{center}
\begin{tabular}{|c|l|}
\hline
$n$ & size of possible dissections \\
\hline
\hline 10 & 4, 7  \\
\hline 12 & 4, 7, 8, 9  \\
\hline 12 & 4, 9, 11  \\
\hline 12 & 4, 11  \\
\hline 12 & 6, 9  \\
\hline 12 & 9  \\
\hline 12 & 11  \\
\hline 13 & 4, 7  \\
\hline 13 & 4, 7, 9, 10  \\
\hline 13 & 4, 7, 10  \\
\hline 13 & 4, 9, 10, 12  \\
\hline 13 & 4, 11, 12  \\
\hline 13 & 9, 10, 12  \\
\hline 13 & 10, 11  \\
\hline
\end{tabular}
\end{center}

A {\em trivial dissection} has triangles of only one size. Apart from
$n = 4$, all trivial dissections are nonseparated.  The following
table lists lower bounds on the number of bitrades that give rise to
the (unique) separated dissection of size $n$ (naturally we allow for
nonseparated solutions to find these trivial dissections).

\begin{center}
\begin{tabular}{|c|l|}
\hline
$n$ & lower bound on number of source bitrades \\
\hline
\hline 4 & 2380591 \\
\hline 8 & 111890 \\
\hline 13 & 1321 \\ 
\hline
\end{tabular}
\end{center}

For each size $n$ we collect examples of dissections with the largest
relative difference in size between the largest and smallest triangle in
the dissection. In all cases the smallest triangle has size $1$ and the
largest triangle  is given in the second column of the table below:

\begin{center}
\begin{tabular}{|c|c|c|c|}
\hline
$n$ & size of largest triangle \\
\hline
\hline 4 & 1 \\
\hline 6 & 2 \\
\hline 7 & 2 \\
\hline 8 & 3 \\
\hline 9 & 4 \\
\hline 10 & 5 \\
\hline 11 & 7 \\
\hline 12 & 9 \\
\hline 13 & 12 \\
\hline 14 & 16 \\
\hline 15 & 21 \\
\hline 16 & 28 \\
\hline 17 & 37 \\
\hline 18 & 49 \\
\hline 19 & 67 \\
\hline 20 & 91 \\
\hline
\end{tabular}
\end{center}

The dissections that give rise to these maximum ratios are shown below
(sorted by dissection size $n$):

\begin{center}
\includegraphics{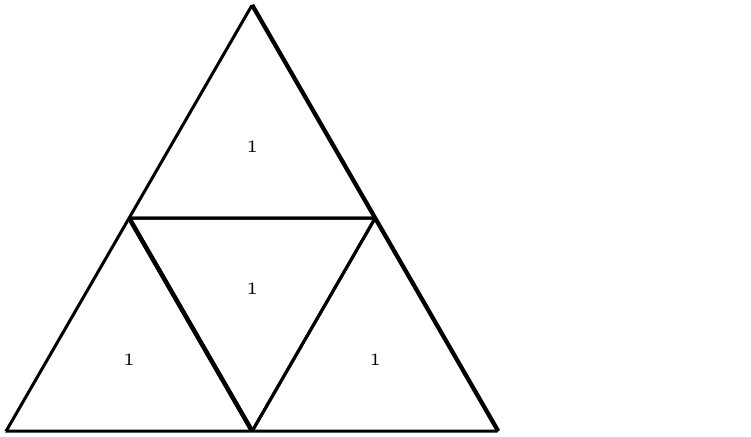}\includegraphics{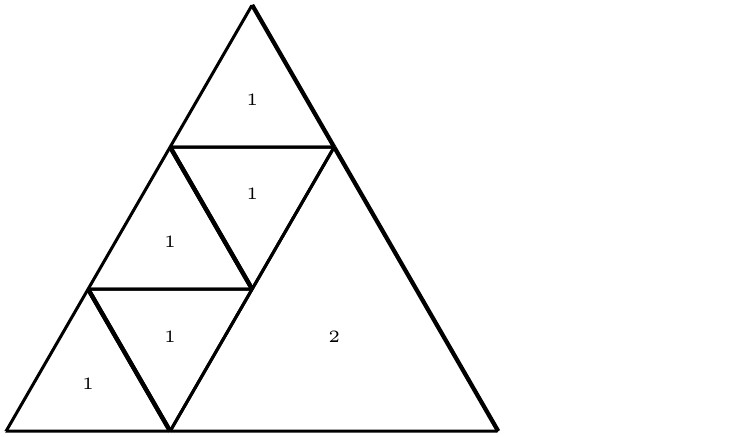}
\end{center}
\begin{center}
\includegraphics{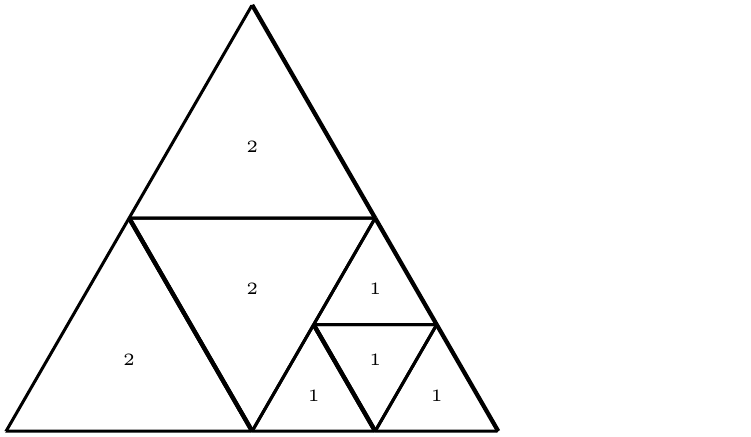}\includegraphics{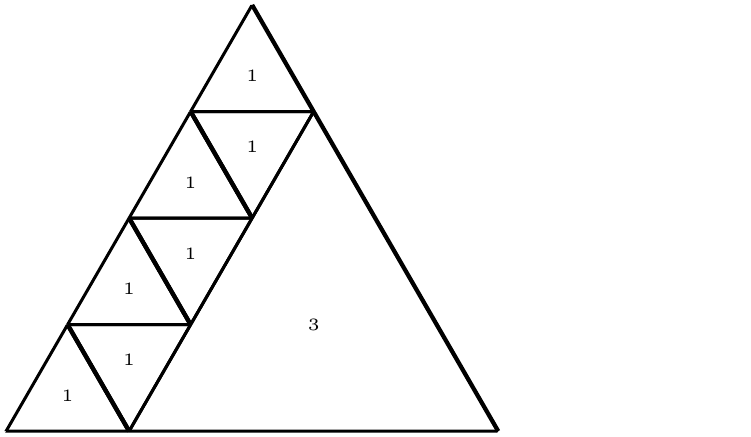}
\end{center}
\begin{center}
\includegraphics{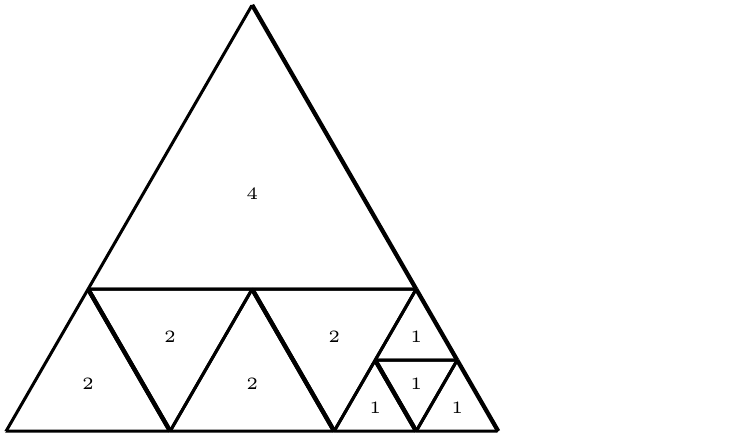}\includegraphics{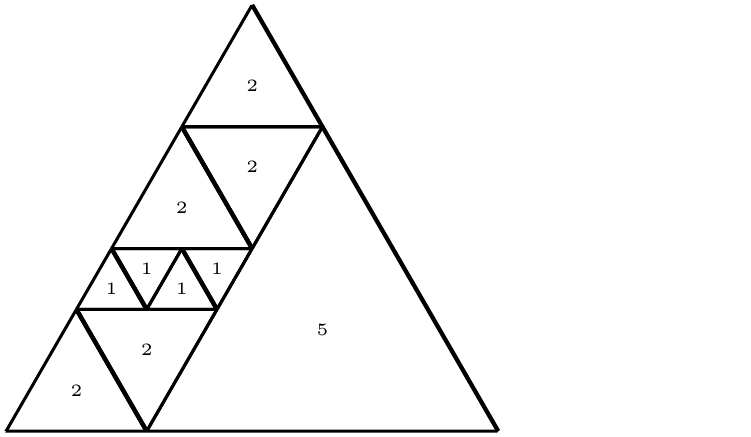}
\end{center}
\begin{center}
\includegraphics{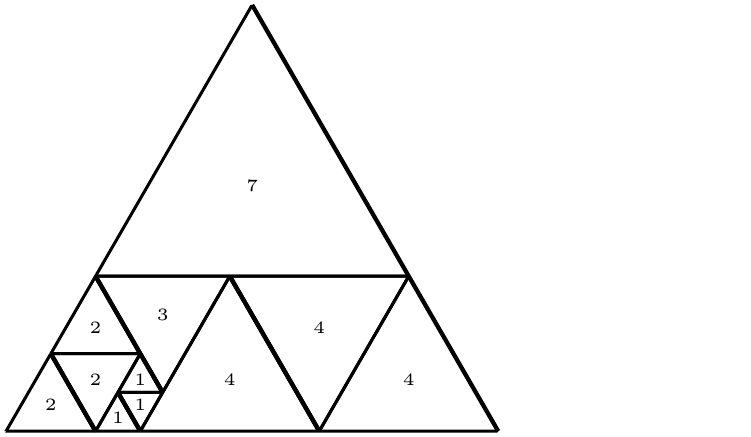}\includegraphics{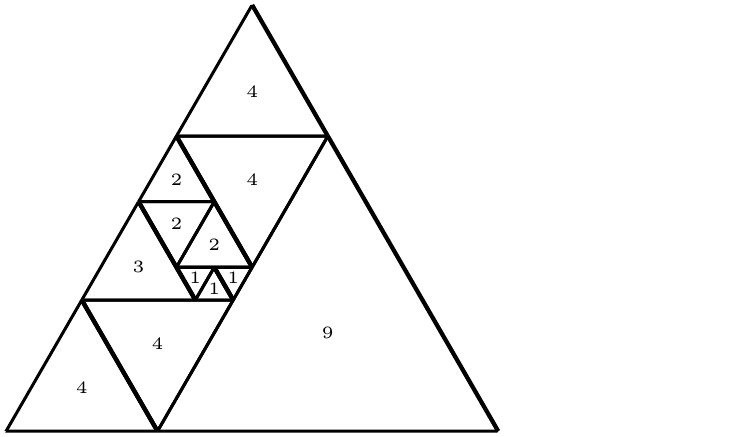}
\end{center}
\begin{center}
\includegraphics{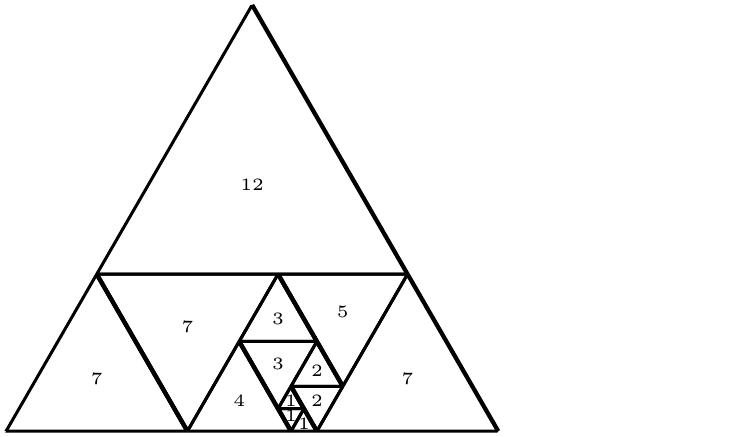}\includegraphics{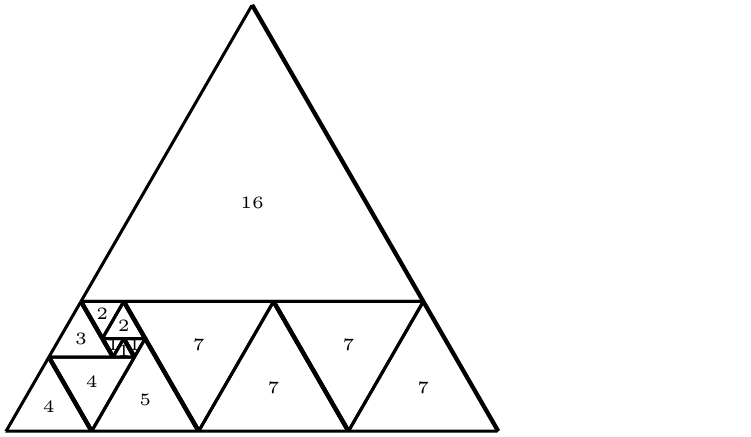}
\end{center}
\begin{center}
\includegraphics{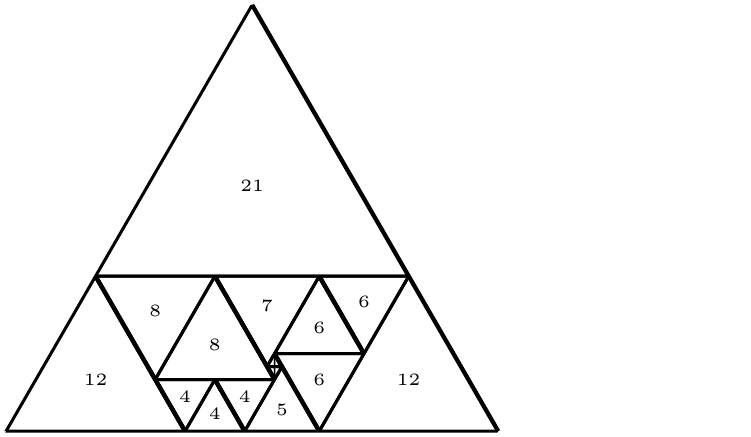}\includegraphics{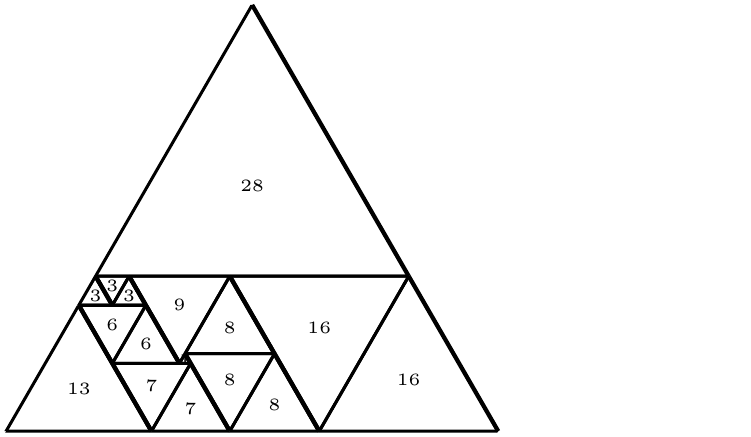}
\end{center}
\begin{center}
\includegraphics{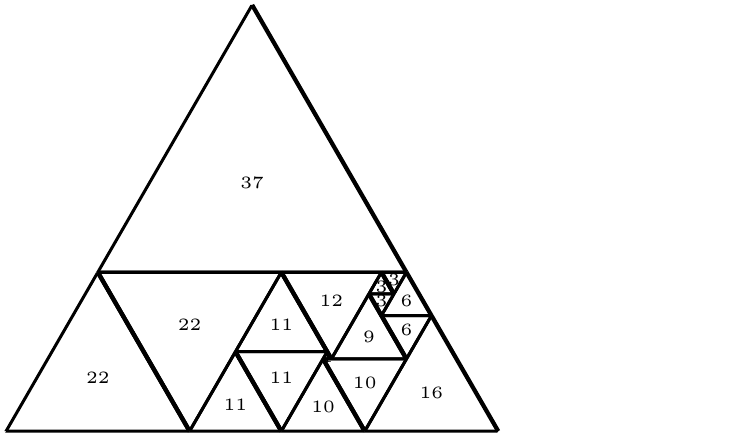}\includegraphics{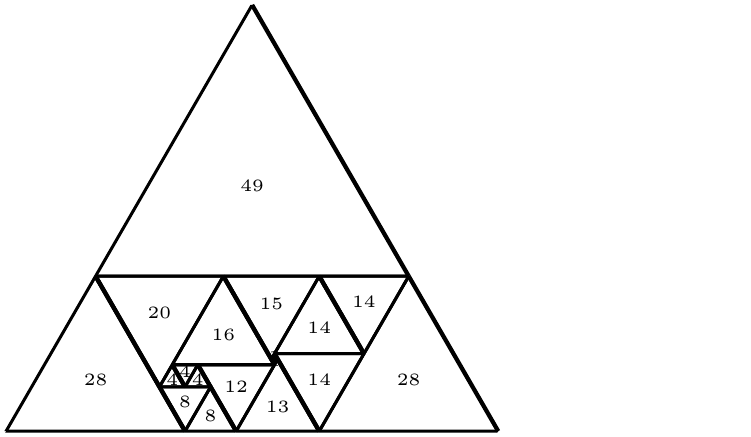}
\end{center}
\begin{center}
\includegraphics{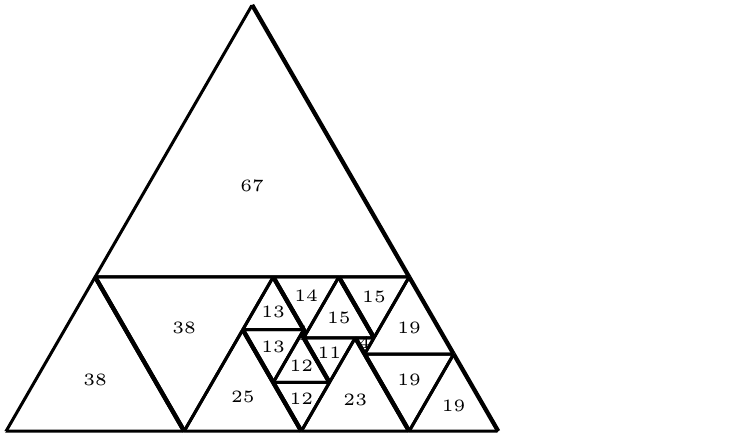}\includegraphics{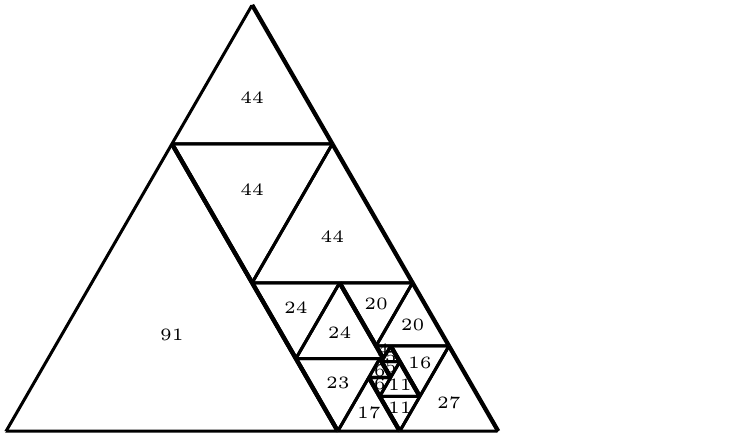}
\end{center}

\clearpage


\clearpage

\section*{Appendix A: Triangle dissections}

Here we present representatives of the isomorphism classes of triangle
dissections of size $n \in \{4, 6, 7, 8, 9, 10 \}$.


\subsection*{$n = 4$}

\includegraphics{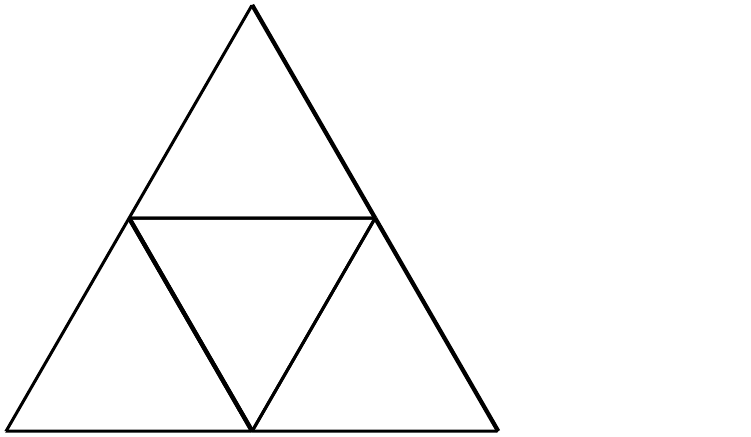}

\subsection*{$n = 6$}

\includegraphics{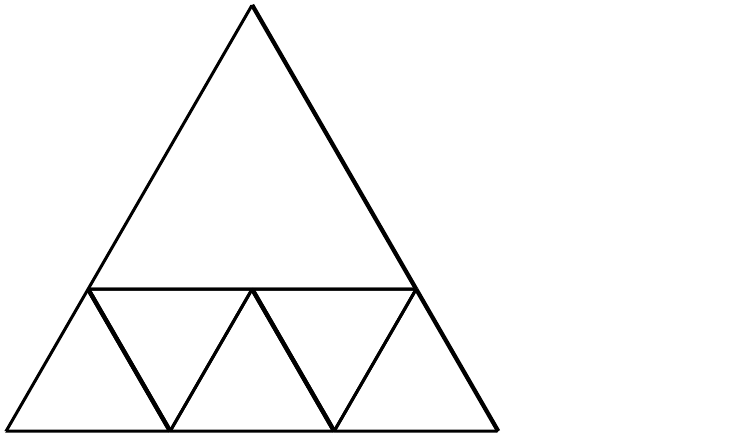}

\subsection*{$n = 7$}

\includegraphics{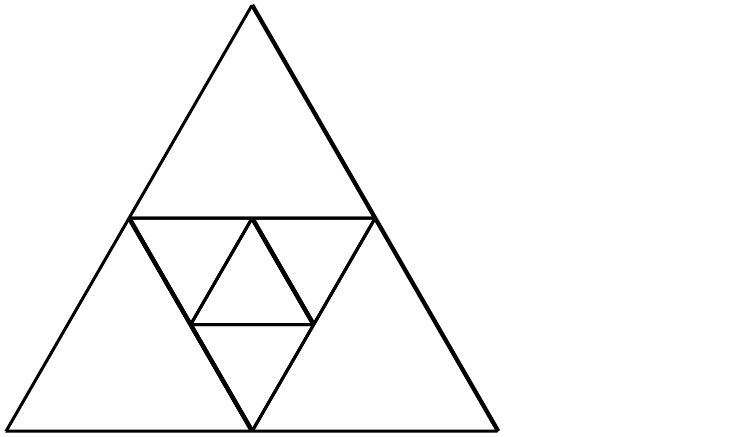}
\includegraphics{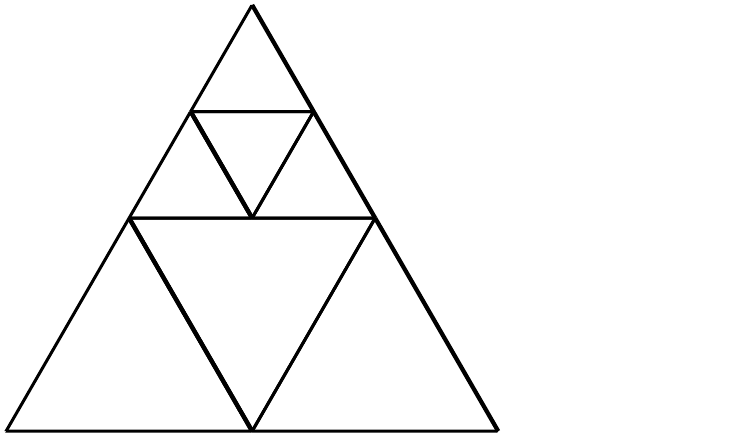}

\subsection*{$n = 8$}

\includegraphics{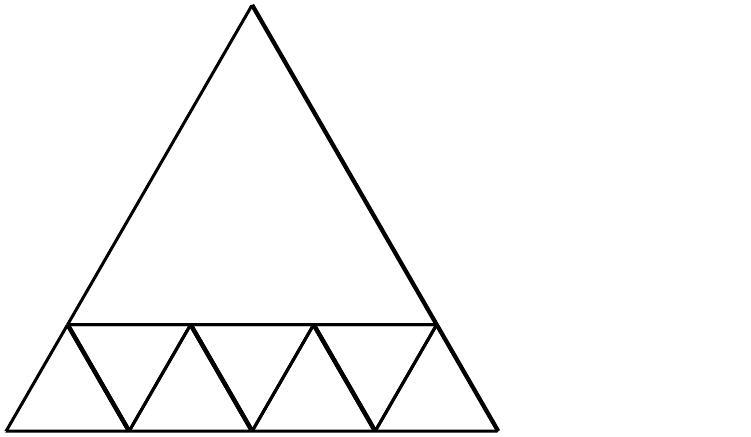}
\includegraphics{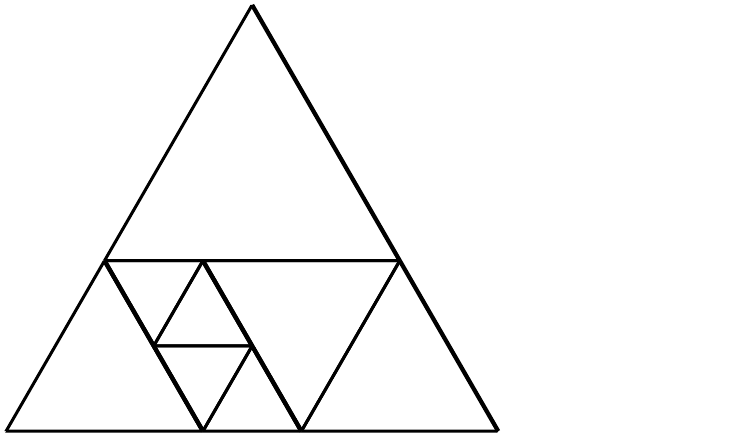}
\includegraphics{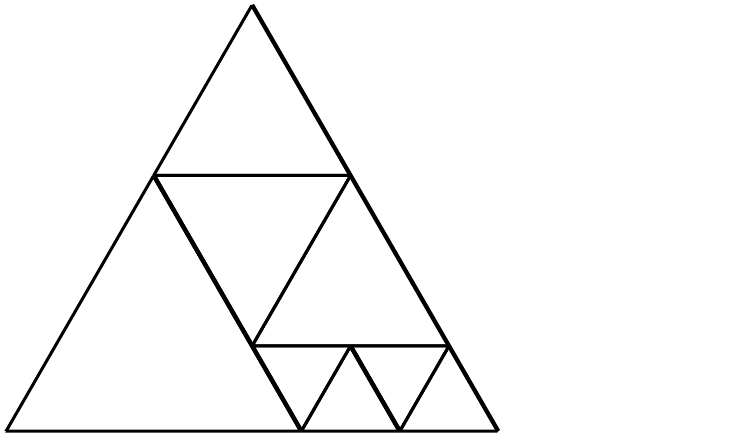}

\subsection*{$n = 9$}

\includegraphics{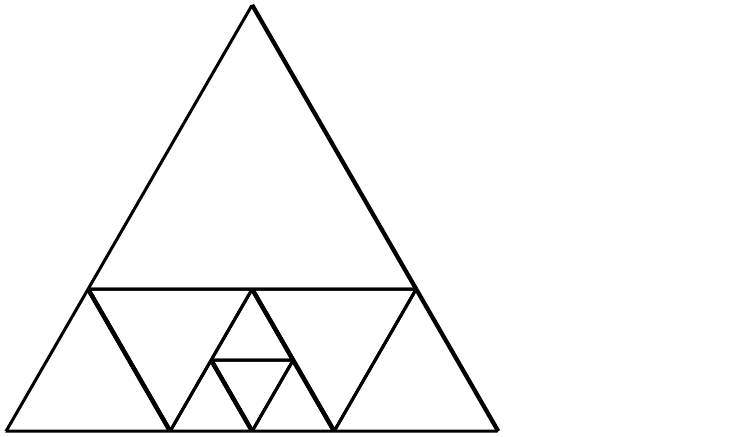}
\includegraphics{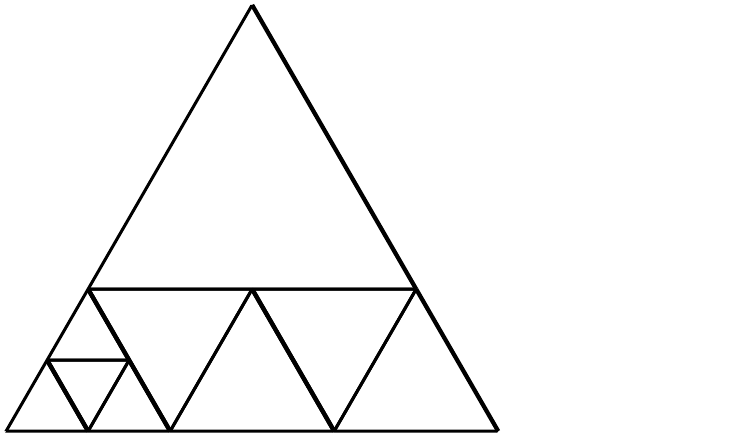}
\includegraphics{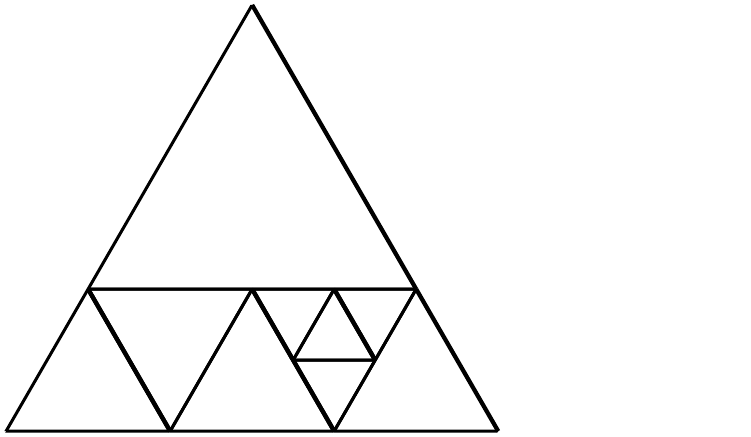}
\includegraphics{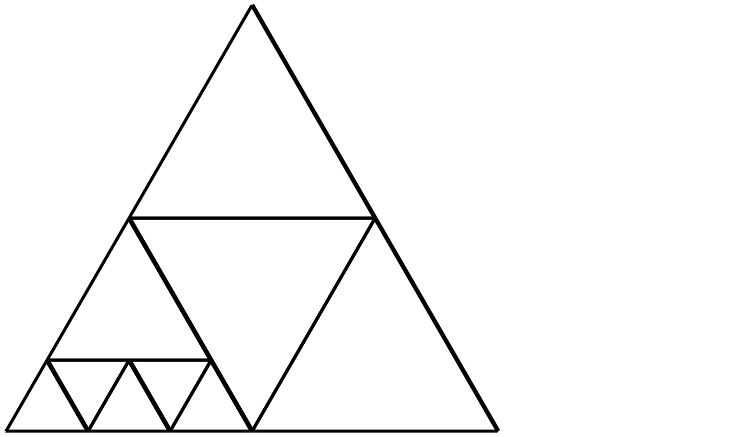}
\includegraphics{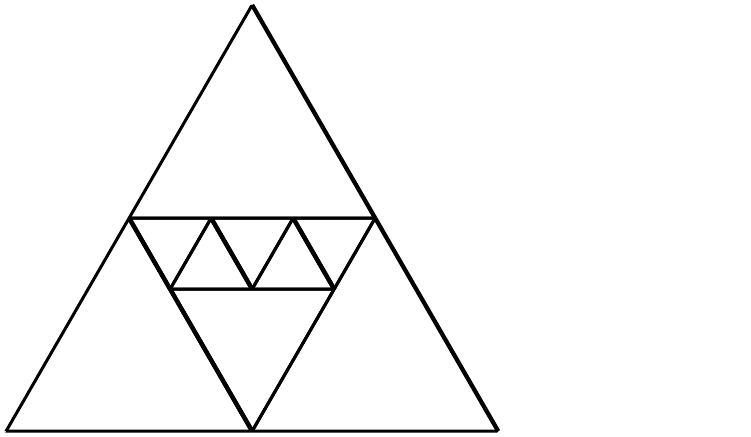}
\includegraphics{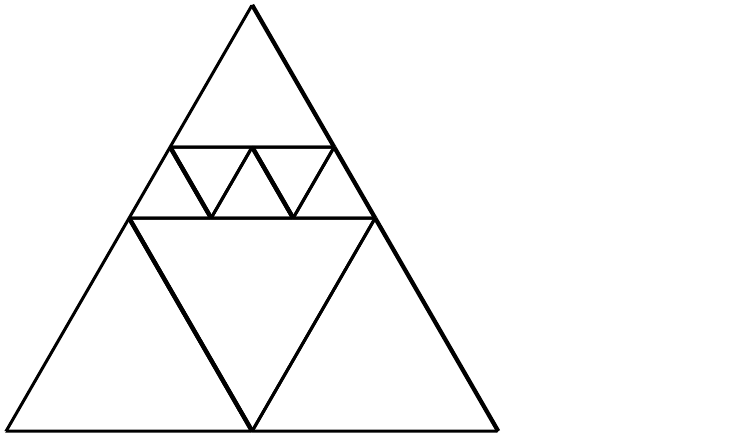}
\includegraphics{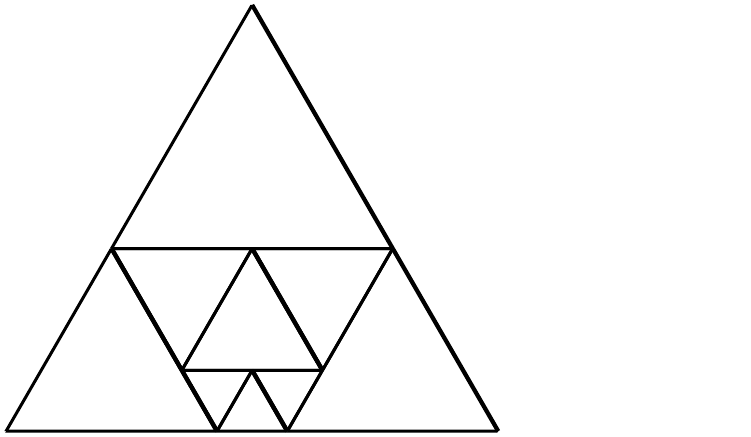}
\includegraphics{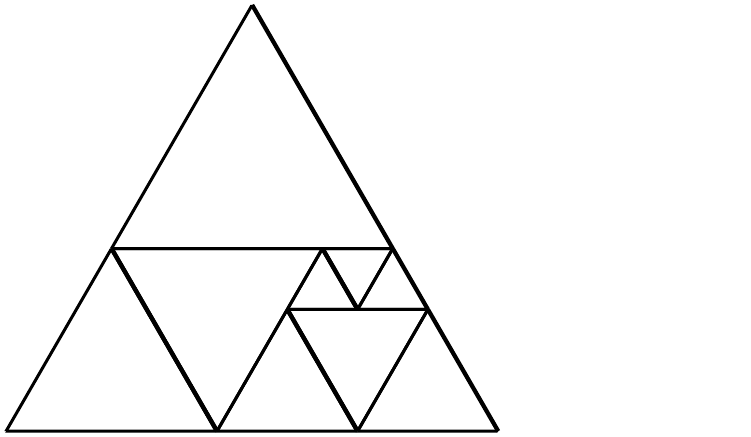}
\includegraphics{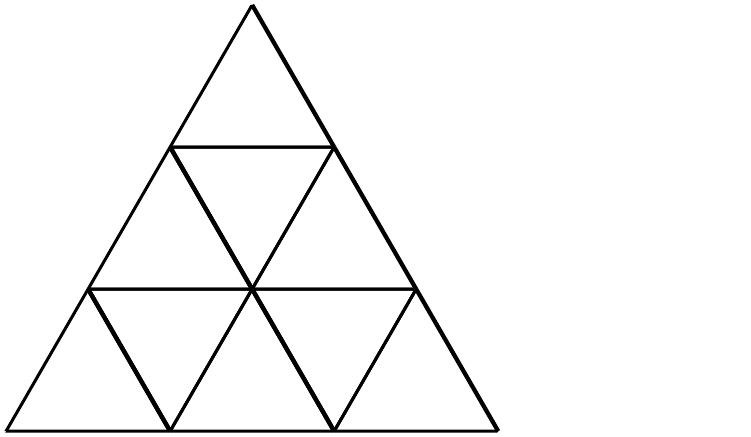}

\subsection*{$n = 10$}

\includegraphics{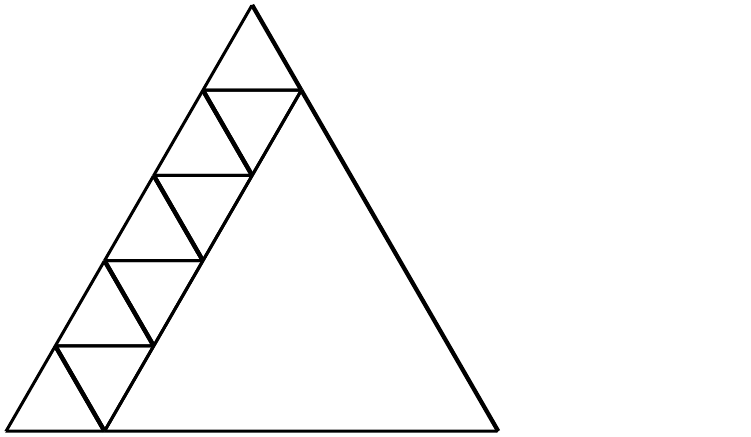}
\includegraphics{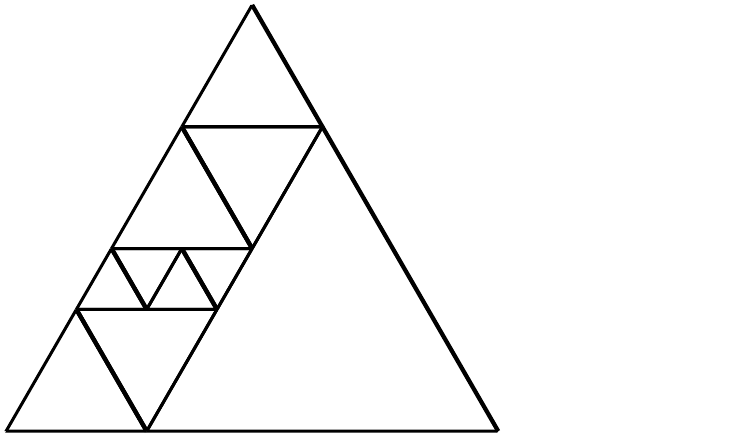}
\includegraphics{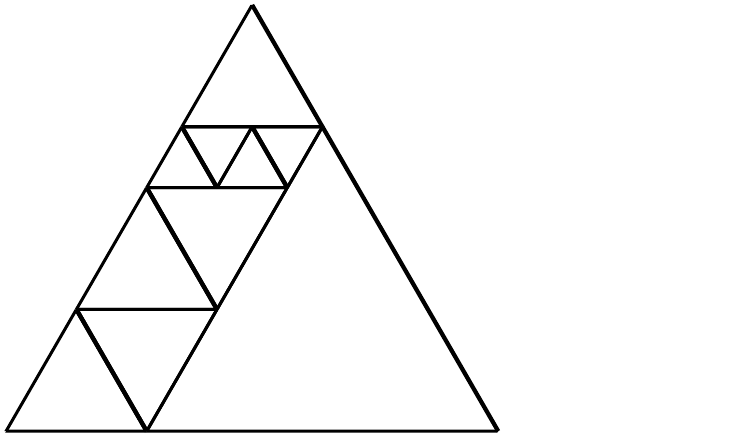}
\includegraphics{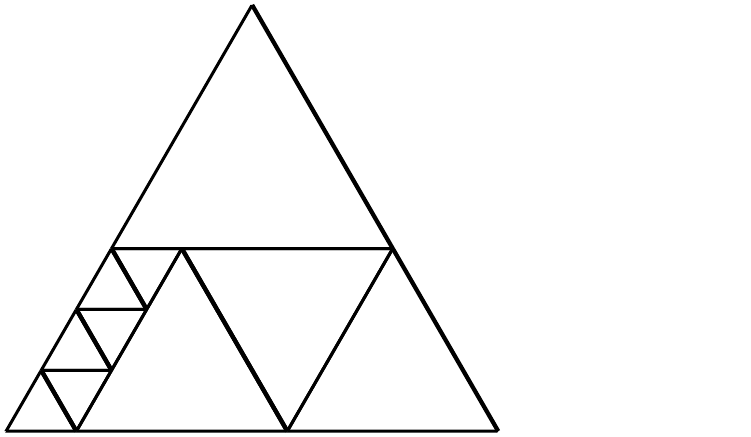}
\includegraphics{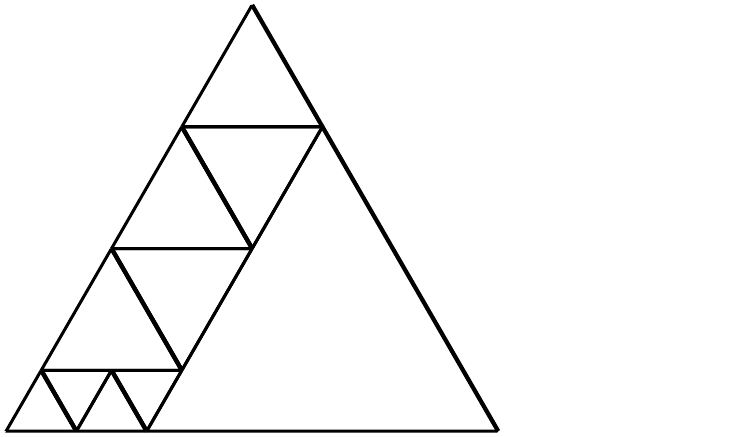}
\includegraphics{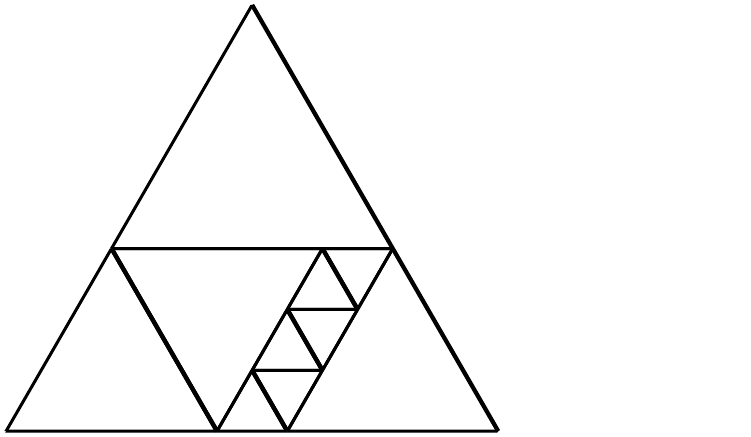}
\includegraphics{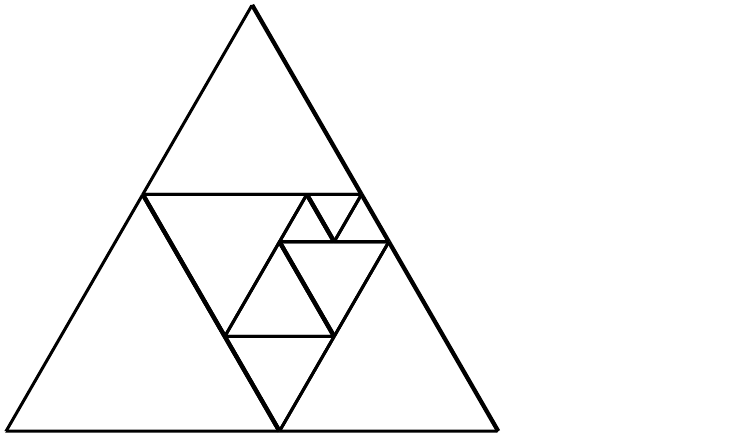}
\includegraphics{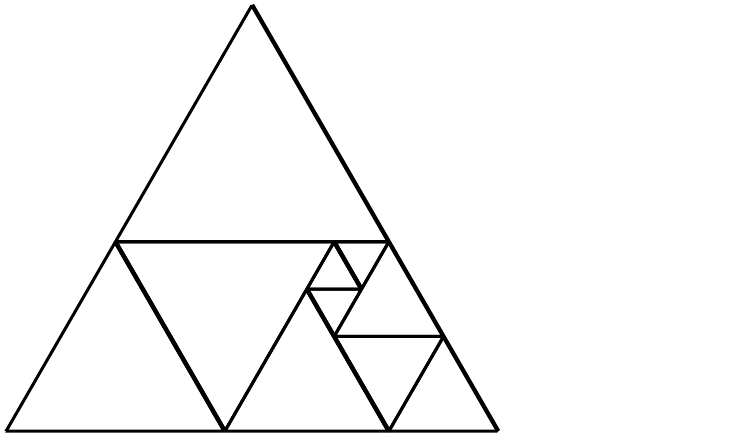}
\includegraphics{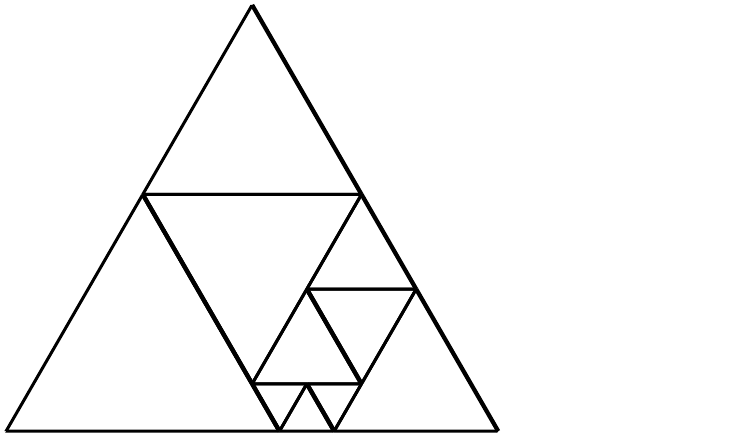}
\includegraphics{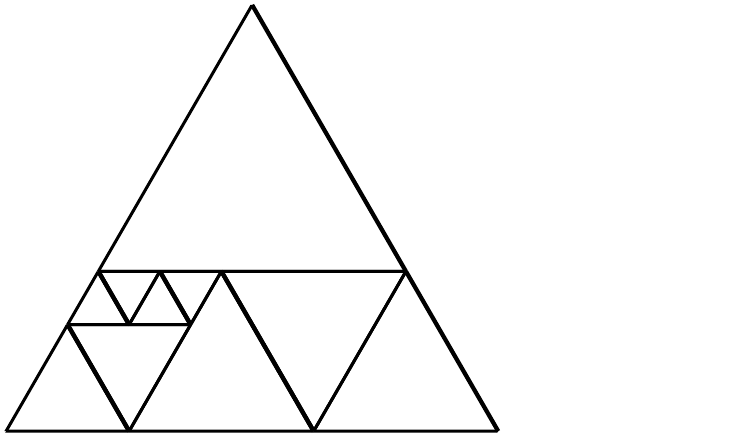}
\includegraphics{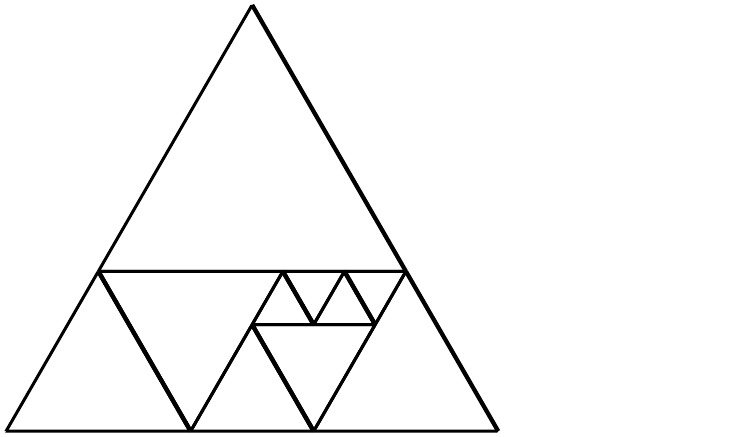}
\includegraphics{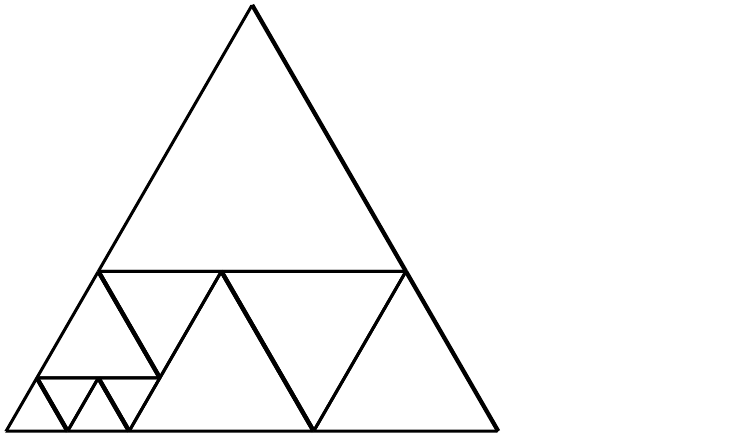}
\includegraphics{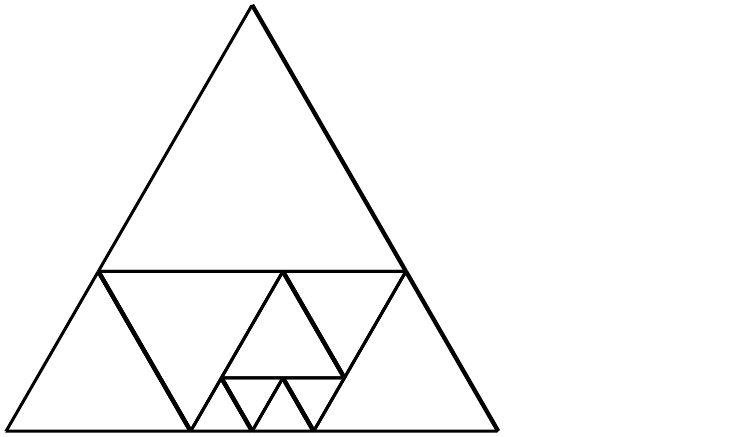}
\includegraphics{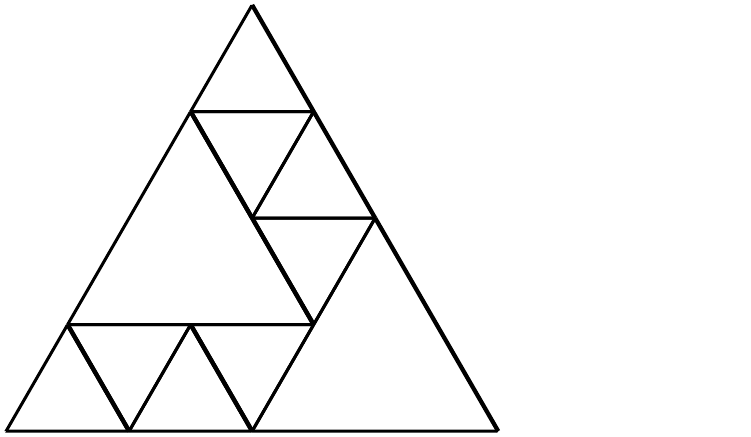}
\includegraphics{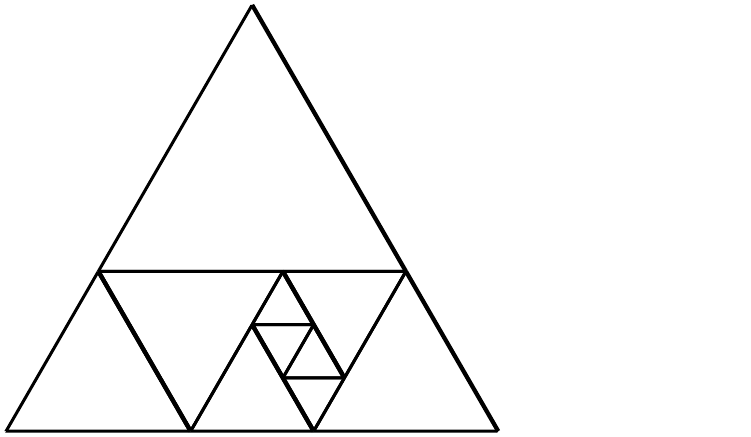}
\includegraphics{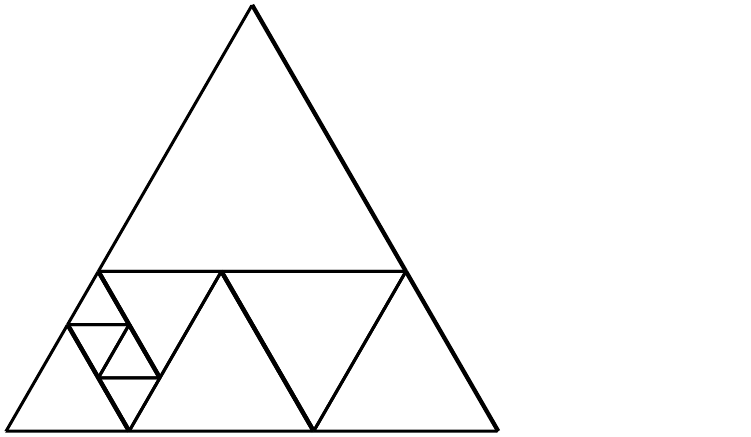}
\includegraphics{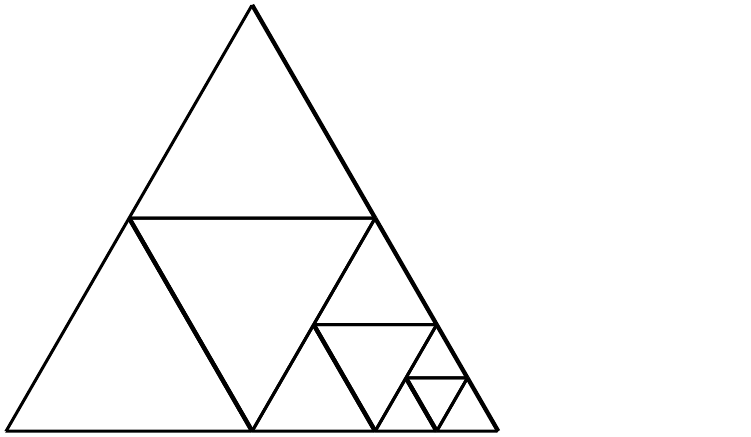}
\includegraphics{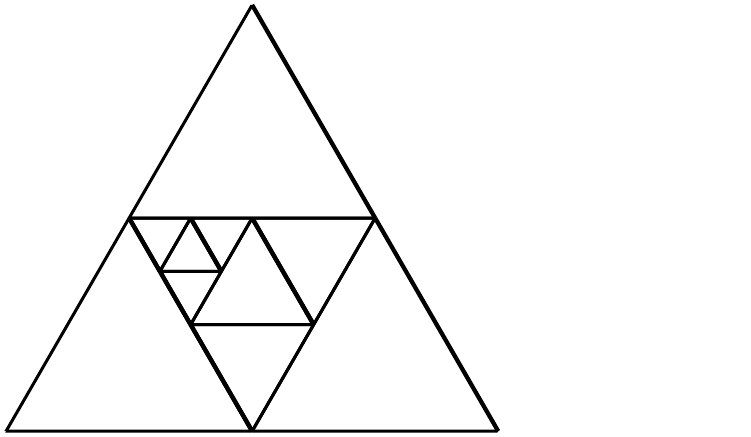}
\includegraphics{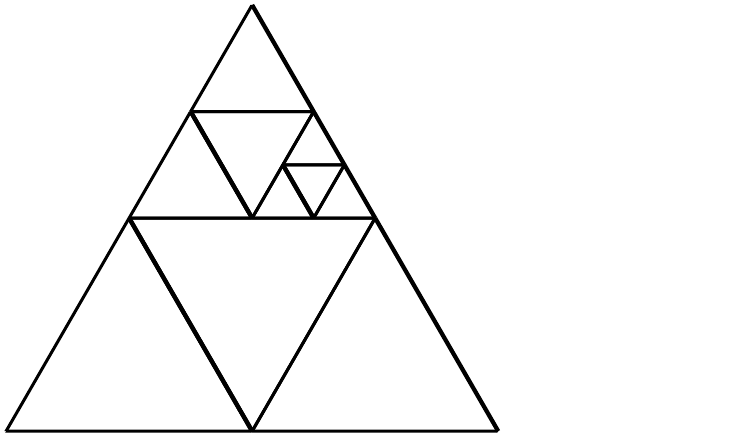}
\includegraphics{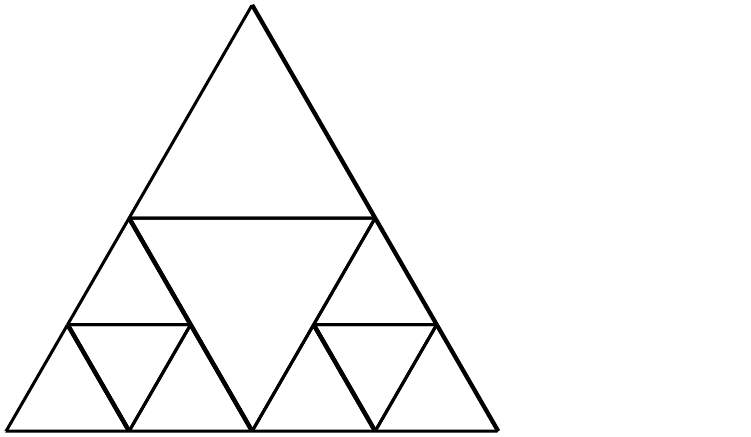}
\includegraphics{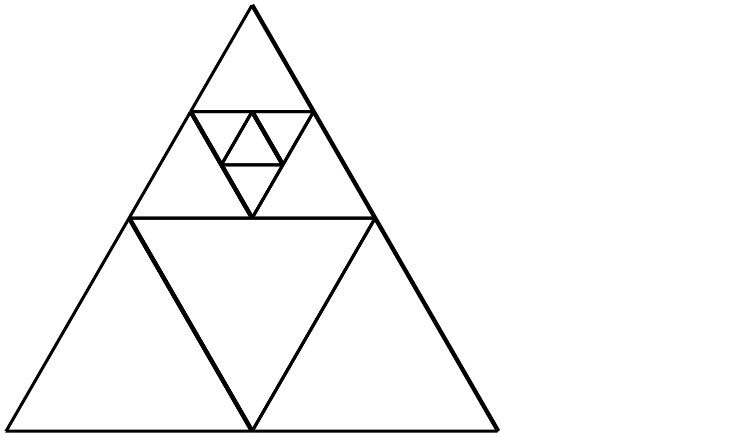}
\includegraphics{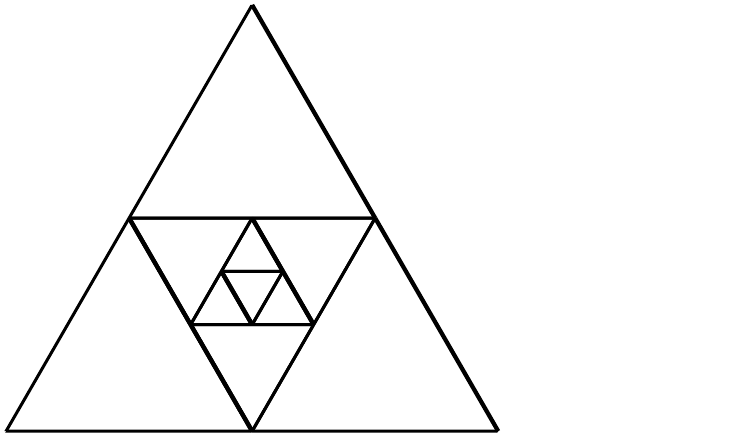}
\includegraphics{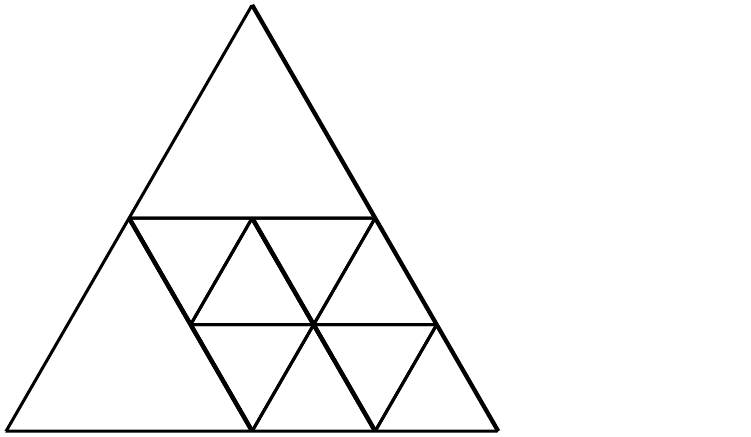}

\end{document}